\documentclass[11pt, leqno]{article}

\usepackage{amssymb}

\newtheorem{thm}{Theorem}[section]
\newtheorem{prop}[thm]{Proposition}
\newtheorem{lemma}[thm]{Lemma}

\newtheorem{Remarks}[thm]{Remarks}
\newtheorem{Remark}[thm]{Remark}

\newtheorem{conjecture}[thm]{Conjecture}
\newtheorem{cor}[thm]{Corollary}

\newcounter{ex}[section]

\newcommand{\Mod}{\mathfrak {Mod}}

\newcommand{\R}{{\bf R}}
\newcommand{\Q}{{\bf Q}}

\newcommand{\Gg}{{\cal G}}

\newcommand{\F}{{\cal F}}

\newcommand{\Spec}{{\rm Spec }\, }
 \renewcommand{\O}{{\cal O}}

\newcommand{\I}{{\cal I}}

\newcommand{\lon}{{\longrightarrow}}
\newcommand{\lo}{{\rightarrow}}

\newcommand{\phil}{f}

\newcommand{\N}{{\cal N}}

\newcommand{\bcf}{\hbox{\boldmath $\cal F$}}

\newcommand{\Fsmile}{\breve F}
\newcommand{\Ksmile}{\breve K}
\newcommand{\Esmile}{\breve E}
\newcommand{\Msmile}{\breve M}

\def\thfill{\null\nobreak\hfill}

\def\endproof{\thfill\vbox{\hrule
  \hbox{\vrule\hbox to 5pt{\vbox to 5pt{\vfil}\hfil}\vrule}\hrule}}

\newenvironment{eq}{\addtocounter{thm}{0}\begin{equation} }{\end{equation}}

\begin{document}
\title{Local models in the ramified case\\ I. The EL-case}
\author{G. Pappas and M. Rapoport}
\date{\ }
\maketitle



\section{Introduction} \label{intro}
\setcounter{equation}{0}


In the arithmetic theory of Shimura varieties it is of interest to have a
model over $\Spec\O_E$, where $E$ is the completion of the reflex field at
some finite prime of residue characteristic $p$. If the Shimura variety is
the moduli space of abelian varieties of PEL type and the level structure
at $p$ is of parahoric type, such a model was defined in [RZ] by posing
the moduli problem over $\Spec\O_E$. In [RZ] it was conjectured that this
model is flat over $\Spec\O_E$, but in [P] it was shown that this is not
always true. It still seems reasonable to expect this flatness property
when the group $G$ defining the Shimura variety splits over an unramified
extension of $\Q_p$, and this is supported by the theorem of G\"ortz [G]
in the EL-case. On the other hand, in [P] it is shown that this conjecture
fails in general if the group $G$ defining the Shimura variety has
localization over $\Q_p$ isomorphic to the group of unitary similitudes
corresponding to a {\it ramified} quadratic extension of $\Q_p$.
Furthermore, a modified moduli problem was proposed in loc. cit. which
defines a closed subscheme of the original model and which stands a better
chance of being flat over $\Spec\O_E$.

Our main purpose of the present paper is to come to grips with the
phenomenon of non-flatness by investigating the simplest case in which it
can occur. As usual, the problem can be reduced to the consideration of
the associated {\it local model} which locally for the \'etale topology
around each point of the special fiber coincides with the model of the
Shimura variety. In the rest of the paper we only consider the local
models which can be defined in terms of linear algebra as schemes over the
spectrum of a complete discrete valuation ring with perfect residue field.
However, in view of the fact that the models proposed in [RZ] are not
``the right ones'' in general we shall term them {\it naive local models}
and reserve the name of {\it local models} for certain closed subschemes
of the naive local models defined in the body of the paper. Both of them
have as generic fibers a closed subscheme of a Grassmannian, and as
special fiber a closed subvariety of the affine partial flag variety over
the residue field corresponding to the fixed parahoric.

As a by-product of our investigations, as they concern the special fibers,
we also obtain several results on the structure of Schubert varieties in
affine Grassmannians and their relation to nilpotent orbit closures.

The simplest case of a naive local  model occurs for the {\it standard
models.} Let us define them. Let $F_0$ be a complete discretely valued
field with ring of integers $\O_{F_0}$ and perfect residue field. Let
$F/F_0$ be a totally ramified extension of degree $e$, with ring of
integers $\O_F$. Let $V$ be a $F$-vector space of dimension $d$, and let
$\Lambda$ be a $\O_F$-lattice in $V$. Choose for each embedding $\varphi$
of $F$ into a separable closure $F_0^{\rm sep}$ of $F_0$ an integer
$r_{\varphi}$ with

\begin{eq}\label{1.1}
0\leq r_{\varphi}\leq d\ \ ,\ \ \forall\varphi\ \ .
\end{eq}
Put $r=\sum_{\varphi}r_{\varphi}$. Then the {\it standard model for $GL_d$
and } ${\bf r}=(r_{\varphi})$, denoted $M(\Lambda, {\bf r})$, parametrizes
the points in the Grassmannian of subspaces ${\cal F}$ of rank $r$ of
$\Lambda$ which are $\O_F$-stable and on which the representation of
$\O_F$ is prescribed in terms of ${\bf r}$ (comp.\ (\ref{2.4})). It is
defined over $\Spec\O_E$ where $E=E(V, {\bf r})$ is the reflex field, the
field of definition of the prescribed representation of $\O_F$. Over
$F_0^{\rm sep}$ the scheme $M=M(\Lambda, {\bf r})$ is isomorphic to the
product over all $\varphi$ of Grassmannians of subspaces of dimension
$r_{\varphi}$ in a vector space of dimension $d$. Over the residue field
$k$ of $\O_E$ the scheme $\overline M=M\otimes_{\O_E}k$ is a closed
subvariety of the affine Grassmannian of $GL_d$ over $k$, and is in fact a
union of Schubert strata which are enumerated by the following dominant
coweights of $GL_d$

\begin{eq}\label{1.2}
{\cal S}^0(r,e,d)= \{\ {\bf s}=s_1\geq \ldots\geq s_d;\ e\geq s_1,\
s_d\geq 0,\ \sum_i {s_i}=r\ \}\ .
\end{eq}
An easy dimension count shows that the dimension of the generic fiber and
the special fiber are unequal, unless all integers $r_{\varphi}$ differ by
at most 1. In this last case we conjecture, and often can prove that $M$
is flat over $\Spec\O_E$. In particular, in the Hilbert-Blumenthal case
($d=2$, and $r_{\varphi}=1$, $\forall\varphi$) the standard model
coincides with the flat model constructed in this case in [DP].  In all
cases when two $r_{\varphi}$ differ from each other by more than one, $M$
is not flat over $\Spec\O_E$.

To analyze the standard model we construct a diagram

\begin{eq}\label{1.3}
M\buildrel\pi\over\longleftarrow \widetilde
M\buildrel\phi\over\longrightarrow N\ \ .
\end{eq}

Here $\widetilde M$ is the $GL_r$-torsor over $M$ which fixes a basis of
the variable subspace ${\cal F}$. The morphism $\phi$ is defined as
follows. Let $\pi$ be a uniformizer of $\O_F$ which satisfies an
Eisenstein polynomial $Q(\pi)=0$. Let

\begin{eq}\label{1.4}
N{=}\{ A\in {\rm Mat}_{r\times r};\ {\rm det}(T\cdot I-A)\equiv
\prod_{\varphi}(T-\varphi(\pi))^{r_{\varphi}}\ ,Q(A){=}0\}  .
\end{eq}
By considering the action of $\pi$ on the variable subspace ${\cal F}$ and
expressing it as a matrix in terms of the fixed basis of ${\cal F}$, we
obtain the morphism $\phi$. We show the following
result (Theorem \ref{phismooth}).

\bigskip\noindent
{\bf Theorem A:} {\it The morphism $\phi$ is smooth of relative dimension
$rd$.}

\bigskip

 Using Theorem A many structure problems on $M$ can be reduced to
corresponding questions on $N$. After a finite extension $\O_E\to \O_K$,
the variety $N$ can be seen as a rank variety in the sense of Eisenbud and
Saltman ([ES]), whereas the special fiber $\overline N=N\otimes_{\O_E}k$
is a subscheme of the nilpotent variety,

\begin{eq}\label{1.5}
\overline N= \{ A\in {\rm Mat}_{r\times r};\ {\rm det}(T\cdot I-A)\equiv
T^r,\ A^{e}=0\}\ .
\end{eq}

Using now results of Mehta and van der Kallen ([M-vdK]) on the structure
of the closures of nilpotent conjugacy classes and basing ourselves on the
methods of Eisenbud and Saltman, we use this reduction procedure to prove
the following result (Theorem \ref{mcan}).

\bigskip\noindent
{\bf Theorem B:} {\it Let $M^{\rm loc}$ be the scheme theoretic closure of
$M\otimes_{\O_E}E$ in $M$. Then

\begin{enumerate}
\item[(i)] $M^{\rm loc}$ is normal and Cohen-Macaulay.
\item[(ii)] The special fiber $\overline M^{\rm loc}$ is reduced, normal
with rational singularities and is the union of all those strata of
$\overline M$ which correspond to ${\bf s}$ in (\ref{1.2}) with ${\bf
s}\leq {\bf r}^{\vee}$ (dual partition to ${\bf r}$).
\item[(iii)] If the scheme $\overline N$ is reduced, then
$M^{\rm loc}=M$ provided that all $r_{\varphi}$ differ by at most 1.
\end{enumerate}
}

We conjecture that the hypothesis made in (iii) is automatically
satisfied, i.e.\ that $\overline N$ is always reduced. This is true by a
classical result of Kostant if $r\leq e$ (in which case the second
condition in (\ref{1.5}) is redundant). For $e<r$ it seems a difficult
problem. In a companion paper to ours, J.~Weyman proves our conjecture for
$e=2$, and for arbitrary $e$ when char $k=0$.

Recall that Lusztig [L] has interpreted certain Schubert varieties in the
affine Grassmannian of $GL_r$ as a compactification of the nilpotent
variety of $GL_r$ (namely the Schubert variety corresponding to the
coweight $(r,0,\ldots, 0)$), compatible with the orbit stratifications of
both varieties. In particular, as used by Lusztig in his paper, all
singularities of nilpotent orbit closures occur in certain Schubert
varieties in the affine Grassmannians.
However, the main thrust of Theorem A (when restricted to the special fibers) goes in the
other direction. Namely, we prove:

\bigskip\noindent
{\bf Theorem C:} {\it Any Schubert variety in the affine Grassmannian of
$GL_d$ is smoothly equivalent to a nilpotent orbit closure for $GL_r$, for
suitable $r$. In particular, it is normal with rational singularities. }

\medskip
We point out the recent preprint by Faltings [F] in which he proves basic
results (like normality) on Schubert varieties in affine flag varieties
for arbitrary reductive groups. \par
 The disadvantage of $M^{\rm loc}$
compared to $M$ is that $M^{\rm loc}$ is not defined by a moduli problem
in general. However, assume that $e=2$ and order the embeddings so that
$r_{\varphi_1}\geq r_{\varphi_2}$. If $r_{\varphi_1}=r_{\varphi_2}$, then
$E=F_0$ and it follows from (iii) above and Weyman's result that $M^{\rm
loc}=M$. If $r_{\varphi_1}>r_{\varphi_2}$, we may use $\varphi_1$ to
identify $E$ with $F$. Then using work of Strickland (see Cor.
\ref{cor5.8}) we can see that $M^{\rm loc}$ is defined inside $M$ by the
following condition on ${\cal F}$,

\begin{eq}\label{1.6}
\wedge^{r_{\varphi_2}+1}(\pi-\varphi_1(\pi)\cdot {\rm Id}\vert {\cal
F})=0\ \ .
\end{eq}

In general, there is a connection to a conjecture of De Concini and
Procesi describing the ideal of the closure of a nilpotent conjugacy
class. Assuming this conjecture to be true and under a technical
hypothesis we can write down a number of conditions which would define
$M^{\rm loc}$ inside $M$. Since these conditions are highly redundant, the
interest of such a description may be somewhat limited, however.

Let $N^{\rm loc}$ be the scheme theoretic closure of $N\otimes_{\O_E}E$ in
$N$. Then the special fiber $\overline N^{\rm loc}$ can be identified with
the closure in $\overline N$ of the nilpotent conjugacy class
corresponding to ${\bf r}^{\vee}$. Let $K$ be the Galois closure of
$F/F_0$ with ring of integers $\O_K$ and residue field $k'$. Then
$N\otimes_{\O_E}\O_K$ has a canonical resolution of singularities,

\begin{eq}\label{1.7}
\mu: {\cal N}\longrightarrow N\otimes_{\O_E}\O_K\ ,
\end{eq}
i.e. a morphism with source a smooth $\O_K$-scheme ${\cal N}$ which is an
isomorphism on the generic fiber and whose special fiber can be identified
with the Springer-Spaltenstein resolution of $\overline N^{\rm loc}
\otimes_kk'$. Using the fact that $\mu\otimes_{\O_K}k'$ is semi-small, one
can calculate the direct image of the constant perverse sheaf on
$\overline{\cal N}$ (Borho, MacPherson [BM], Braverman, Gaitsgory [BG]).
Transporting back the result to $M$ via Theorem A, we obtain the following
result (comp.\ Theorem \ref{theorem6.1}).

\bigskip\noindent
{\bf Theorem D:} {\it Let $R\psi'$ denote the complex of nearby cycles of
the $\O_K$-scheme $M\otimes_{\O_E}\O_K$. Then there is the following
identity of perverse sheaves pure of weight 0 on $\overline M\otimes_kk'$,

$$ R\psi'[{\rm dim}\ \overline M]\hbox{$({1\over 2}{\rm dim}\ \overline
M)$}=\bigoplus_{{\bf s}\leq {\bf r}^{\vee}} K_{{\bf r}^{\vee},{\bf
s}}\cdot IC_{M_{\bf s}\otimes_kk'}\ \ . $$ Here $K_{{\bf r}^{\vee},{\bf
s}}$ is a Kostka number. }

\medskip
We refer to section \ref{section6} for an explanation of the notation and
for the question of descending this result from ${\cal O}_K$ to ${\cal
O}_E$.

Via Theorem A, the resolution of singularities (\ref{1.7}) is closely
related to an analogous resolution of singularities of $M\otimes_{{\cal
O}_E} {\cal O}_K$ whose special fiber can be identified with the Demazure
resolution of the affine Schubert variety corresponding to the coweight
${\bf r}^{\vee}$ of $GL_d$. As pointed out to us by Ng\^o, from this
identification one obtains another formula for the complex of nearby
cycles $R\psi'$.

This concludes the discussion of our results on the standard models for
$GL_d$ and ${\bf r}$. The general naive local models of $EL$-type are
projective schemes $M^{\rm naive}$ over $\Spec\O_E$ which at least over
$\Spec\O_{\breve E}$ (where $\breve E$ is the completion of the maximal
unramified extension of $E$) are closed subschemes of products of standard
models. We define closed subschemes $M^{\rm loc}$ of $M^{\rm naive}$ by
demanding that the projection into any standard model lies in the flat
closure considered above. The modification of the flatness conjecture of
[RZ] in the present case is that $M^{\rm loc}$ is flat over $\Spec\O_E$.
We have nothing to say about this conjecture, except to point out that the
description of $M^{\rm loc}(k)$ implicit in [KR], in terms of
$\mu$-permissible elements in the Iwahori-Weyl group is correct, as
follows from Theorem B, (ii).

\bigskip

In the original version of this paper Theorem A was proved by exhibiting
affine charts around the worst singularities of $M$ over which the
morphisms $\pi$ and $\phi$ can be made completely explicit. The present
simple proof of Theorem A is based on ideas from [FGKV]. We thank D.\
Gaitsgory for pointing it out to us. We would also like to thank
G.\,Laumon, T.\,Haines and B.C.\,Ng\^o for helpful discussions on the
material of sections \ref{section5A} and \ref{section6}. We are grateful
to J.~Weyman for his interest in our conjecture and for making his results
available in a companion paper to ours. We thank G.~Pfister for computer
calculations using the symbolic algebra package REDUCE, which gave us the
courage to elevate an initially naive question to the rank of a
conjecture. We also thank the Max-Planck-Institut Bonn for its hospitality
and support. The first named author was also partially supported by NSF
grant DMS99-70378 and by a Sloan Research Fellowship.

\bigskip\noindent
{\it General notational convention:} If $X$ is a scheme over Spec $R$ and
$R'$ is an $R$-algebra, we often write $X\otimes_RR'$ or $X_{R'}$ for
$X\times_{\Spec R}\Spec R'$.

\section{Standard models for $GL_d$} \label{GLd}
\setcounter{equation}{0}

In this section and the sections \ref{grassmannian} -- \ref{section6} we
will use the following notation. Let $F_0$ be a complete discretely valued
field with ring of integers ${\cal O}_{F_0}$ and uniformizer $\pi_0$, and
perfect residue field. Let $F$ be a totally ramified separable extension
of degree $e$ of $F_0$, with ring of integers ${\cal O}_F$. Let $\pi$ be a
uniformizer of ${\cal O}_F$ which is a root of the Eisenstein polynomial
\begin{eq}\label{2.1}
Q(T)=T^e+\sum_{k=0}^{e-1}b_kT^k,\ \ b_0\in
\pi_0\cdot \O^{\times}_{F_0},\ b_k\in
(\pi_0).\end{eq}

Let $V$ be an $F$-vector space of dimension $d$ and $\Lambda$ an
$\O_F$-lattice in $V$. We fix a separable closure $F_0^{\rm sep}$ of
$F_0$. Finally, we choose for each embedding $\varphi:F\to F_0^{\rm sep}$
an integer $r_{\varphi}$ with

\begin{eq}\label{2.2}
0\leq r_{\varphi}\leq d\ \ .
\end{eq}
Associated to these data we have the {\it reflex field} $E$, a
finite extension of $F_0$ contained in $F_0^{\rm sep}$ with

\begin{eq}\label{2.3}
{\rm Gal}(F_0^{\rm sep}/E) =\{ \sigma\in {\rm Gal}(F_0^{\rm
sep}/F_0);\ r_{\sigma\varphi}=r_{\varphi},\ \forall \varphi\}\ \ .
\end{eq}
Let $\O_E$ be the ring of integers in $E$. We now formulate a
moduli problem on $({\rm Sch} / \Spec\O_E)$:

\begin{eqnarray}\label{2.4}
 M(S) &=&\{\
 {\F}\subset\Lambda\otimes_{\O_{F_0}}{\O}_S;\hbox{ a
 ${\O}_F\otimes_{\O_{F_0}} {\O}_S$-submodule,}\\  &&
\hbox{which is locally on
 $S$ a direct summand as}\nonumber\\ &&\hbox{ ${\O}_S$-module,
 with ${\rm{det}}(a\ |\ {\cal F})= \prod\limits_{\varphi}
 \varphi(a)^{r_{\varphi}}\}.$}\nonumber
\end{eqnarray}
\medskip

The last identity
 is meant as an identity of polynomial functions on ${\O}_F$
 (comp.\ [K], [RZ]).

It is obvious that this functor is representable by a projective scheme
$M=M(\Lambda, {\bf r})$ over $\Spec{\cal O}_E$. This scheme is called the
{\it standard model for} $GL_d$ {\it corresponding to} ${\bf
r}=(r_{\varphi})_{\varphi}$ {\it (and to} $F_0, F$ {\it and} $\pi$).

Let us analyze the geometric general fiber and the special fiber of $M$.
We have a decomposition
\begin{eq}\label{2.5}
F\otimes_{F_0}F_0^{\rm sep} = \bigoplus_{\varphi:F\to F_0^{\rm
sep}} F_0^{\rm sep}\ \ .
\end{eq}
Correspondingly we get a decomposition of $V\otimes_{F_0} F_0^{\rm
sep}$ into $F_0^{\rm sep}$-vector spaces

\begin{eq}\label{2.6}
V\otimes_{F_0} F_0^{\rm sep} =\bigoplus_{\varphi} V_{\varphi}\ \ .
\end{eq}
Each summand is of dimension $d$. The determinant condition in (\ref{2.4})
can now be interpreted as saying that $M\otimes_{\O_E}F_0^{\rm sep}$
parametrizes subspaces ${\cal F}_{\varphi}$ of $V_{\varphi}$, one for each
$\varphi$, of dimension $r_{\varphi}$. In other words,
\begin{eq}\label{2.7}
M\otimes_{{\O}_E} F_0^{\rm sep}=\prod\limits_{\varphi}\
 {\rm{Grass}}_{r_{\varphi}}(V_{\varphi})\ \ .
\end{eq}
In particular,

\begin{eq}\label{2.8}
\dim\, M\otimes_{\O_E}E= \sum_\varphi\ r_\varphi(d-r_\varphi)\ .
\end{eq}
Denote by $k$ the residue field of $\O_E$. Let us consider $\overline
M=M\otimes_{{\O}_E}k$. Put

\begin{eq}\label{2.9}
W=\Lambda\otimes_{\O_{F_0}}k\ \ ,\ \ \Pi=\pi\otimes
 {\rm{id}}_{k}\ \ . \end{eq}

Then $W$ is a $k$-vector space of dimension $de$ and $\Pi$ is a nilpotent
endomorphism with $\Pi^e=0$. The conditions that the subspace ${\cal
F}\subset W$ gives a point of $\overline{M}$ translate into the following:
${\cal F}$ is $\Pi$-stable, $\dim_k {\cal F}=r:=\sum_\varphi r_\varphi$
and $\det(T-\Pi\, \vert\, {\cal F})\equiv T^r$. In other words
$M\otimes_{{\O}_E}k$ is the closed subscheme of $\Pi$-stable subspaces
${\cal F}$ in ${\rm{Grass}}_r(W)$ which satisfy the above condition on the
characteristic polynomial. We point out that the $k$-scheme
$\overline{M}=M\otimes_{{\O}_E}k$ only depends on $r$, not on the
partition $(r_{\varphi})$ of $r$.

\medskip

\section{Relation to the affine Grassmannian}\label{grassmannian}
\setcounter{equation}{0}

We denote by $\widetilde{{\rm Grass}}_{k }$ the affine Grassmannian over
$k$ associated to $GL_d$. Recall ([BL]) that this is the Ind-scheme over
$\Spec k$ whose $k$-rational points parametrize the $k \lbrack\lbrack
\Pi\rbrack\rbrack$-lattices in $k ((\Pi))^d$. Here $k
\lbrack\lbrack\Pi\rbrack\rbrack$ denotes the power series ring in the
indeterminate $\Pi$ over $k $. On $\widetilde{\rm Grass}_{k }$ we have an
action of the group scheme $\tilde{\cal G}$ over $k $ with $k $-rational
points equal to $GL_d(k \lbrack\lbrack\Pi\rbrack\rbrack)$. The orbits of
this action are finite-dimensional irreducible locally closed subvarieties
(with the reduced scheme structure) which are parametrized by the dominant
coweights of $GL_d$, i.e.\ by $d$-tuples of integers ${\bf s}=(s_1,
\ldots, s_d)$ with
\begin{eq}\label{3.1}
s_1\geq \ldots \geq s_d\ \ .
\end{eq}
Furthermore ([BL]), if $\O_{\bf s}$ denotes the orbit corresponding to
${\bf s}$, we have
\begin{eq}\label{3.2}
{\rm dim}\ \O_{\bf s} =\langle {\bf s}, 2\varrho\rangle\ \ \hbox{and}\ \
\O_{\bf s'}\subset\ {\rm closure}\ (\O_{\bf s})\Leftrightarrow {\bf
s'}\leq {\bf s}\ \ .
\end{eq}
Here $2\varrho =(d-1, d-3, \ldots, 1-d)$ and $\langle\ ,\ \rangle$ denotes
the standard scalar product on $\R^d$. Furthermore, ${\bf s'}\leq{\bf s}$
denotes the usual partial order on dominant coweights, i.e.\ $${\bf
s'}\leq {\bf s}\Leftrightarrow s'_1\leq s_1,\ s'_1+s'_2\leq
s_1+s_2,\ldots, s'_1+\ldots + s'_d=s_1+\ldots +s_d\ \ .$$ Let us fix an
isomorphism of $k [[\Pi]]$-modules
\begin{eq}\label{3.3}
\Lambda\otimes_{\O_{F_0}}k  \simeq  (k [[\Pi]] /\Pi^e)^d\ \ .
\end{eq}
Then we obtain a closed embedding (of Ind-schemes),
\begin{eq}\label{3.4}
\iota : M\otimes_{\O_E}k  \longrightarrow \widetilde{{\rm Grass}}_{k }\ \
.
\end{eq}
On $k $-rational points, $\iota$ sends a point of $M(k )$, corresponding
to a $k [\Pi]$-submodule ${\cal F}$ of $\Lambda\otimes_{\O_{F_0}}k  =(k
[[\Pi]] /\Pi^e)^d$, to its inverse image $\tilde{\cal F}$ in $k
[[\Pi]]^d$, which is a $k [[\Pi]]$-lattice contained in $k [[\Pi]]^d$,

\begin{eq}\label{3.5}
\matrix{
\tilde{\cal F}&\subset & k [[\Pi]]^d\cr
\big\downarrow&&\big\downarrow\cr
{\cal F}&\subset & (k [[\Pi]]/\Pi^e)^d&.
\cr}
\end{eq}
The embedding $\iota$ is equivariant for the action of $\tilde{\cal G}$ in
the following sense. Consider the smooth group scheme ${\cal G}$ over
$\Spec\O_{F_0}$,
\begin{eq}\label{3.6}
{\cal G} =\underline{\rm Aut}_{\O_F}(\Lambda)\ \ .
\end{eq}
In fact, we will only need the base change of ${\cal G}$ to $\Spec\O_E$
which we denote by the same symbol. The group scheme ${\cal G}$ acts on
$M$ by
\begin{eq}\label{3.7}
(g, {\cal F})\longmapsto g({\cal F})\ \ .
\end{eq}
Let
\begin{eqnarray}\label{3.8}
\overline{\cal G}
= {\cal G}\otimes_{\O_E}k
&=&
\underline{\rm Aut}_{k[\Pi]/\Pi^e} (\Lambda\otimes_{\O_{F_0}}k )\\
&\simeq &GL_d(k[[\Pi]]/\Pi^e)\ \ .\nonumber
\end{eqnarray}
In this way $\overline{\cal G}$ becomes a factor group of $\tilde {\cal
G}$ and the equivariance of $\iota$ means that the action of $\tilde{\cal
G}$ stabilizes the image of $\iota$, that the action on this image factors
through $\overline{\cal G}$ and that $\iota$ is $\overline{\cal
G}$-equivariant.

A point of $M$ with values in a field extension $k'$ of $k$, corresponding
to a $\Pi$-stable subspace ${\cal F}$ of $\Lambda\otimes_{\O_{F_0}}k'$,
has image in $\O_{{\bf s}_{\cal F}}$, where ${\bf s}_{\cal F}$ is the
Jordan type of the nilpotent endomorphism $\Pi\mid {\cal F}$. It follows
that the orbit decomposition of $\overline M=M\otimes_{{\cal O}_E}k$ under
the action of $\overline{\cal G}$ has the form
\begin{eq}\label{3.9}
\overline M =\bigcup_{\bf s} M_{\bf s}\ \ .
\end{eq}
Here $M_{\bf s} \simeq \O_{\bf s}$ via $\iota$ and ${\bf s}$ ranges over
the subset of (\ref{3.1}) given by
\begin{eq}\label{3.10}
{\cal S}^0(r,e,d){=}\{{\bf s}{ =}s_1\geq s_2\geq \ldots \geq s_d;\
 e\geq s_1,\  s_d\geq 0,\ \sum_i s_i{=}r\}\  ,
\end{eq}
i.e.\ the partitions of $r$ into at most $d$ parts bounded by $e$.
In all of the above we have ignored nilpotent elements.

\begin{prop}\label{prop3.1}
The special fibre $\overline M=M\otimes_{{\cal O}_E}k$ is irreducible of
dimension $dr-ec^2-(2c+1)f$.
\end{prop}
Here we have written $r=c\cdot e+f$ with $0\leq f<e$.

\medskip\noindent
\begin{Proof}
Among the coweights in ${\cal S}^0(r,e,d)$ there is a unique maximal one,
\begin{eq}\label{3.11}
{\bf s}_{\rm max}={\bf s}_{\rm max} (r,e)= (e,\ldots, e,f,0\ldots 0)=
(e^c,f)\ \ .
\end{eq}
Hence $M_{{\bf s}_{\rm max}}$ is open and dense in $\overline M$. Its
dimension is equal to $\langle {\bf s}_{\rm max}, 2\varrho\rangle$, which
gives the result.
\end{Proof}\endproof

\bigskip Sometimes for convenience we shall number the embeddings $\varphi$
in such a way that the $r_i=r_{\varphi_i}$ form a decreasing sequence
${\bf r}=(r_1\geq r_2\geq \ldots \geq r_e)$. Then ${\bf r}$ is a partition
of $r$ into at most $e$ parts bounded by $d$.

Let ${\bf r}_{\rm min}={\bf r}_{\rm min}(r,e)={\bf s}^{\vee}_{\rm max}$ be
the dual partition to ${\bf s}_{\rm max}$, i.e.
\begin{eq}\label{3.12}
{\bf r}_{\rm min}=(c+1,\ldots, c+1, c,\ldots, c)= ((c+1)^f, c^{e-f})\ \ .
\end{eq}

\begin{prop}\label{prop3.2}
We have $${\rm dim}\ M(\Lambda, {\bf r})\otimes_{\O_E}E\leq {\rm dim}\
M(\Lambda, {\bf r})\otimes_{\O_E}k\ \ ,$$ with equality if and only if
${\bf r}= {\bf r}_{\rm min}(r,e)$ (after renumbering ${\bf r}$), i.e.\ iff
all $r_{\varphi}$ differ by at most one.
\end{prop}

\begin{Proof}
If ${\bf r}={\bf r}_{\rm min}(r,e)$, then
\begin{eqnarray*}
{\rm dim}\ M(\Lambda, {\bf r})\otimes_{\O_E}E &=& dr-\sum_\varphi
r^2_{\varphi}=dr-f(c+1)^2-(e-f)c^2\\ &=& dr-f((c+1)^2-c^2)-ec^2\\ &=&
dr-ec^2-(2c+1)f\\ &=& {\rm dim}\ M(\Lambda, {\bf r}) \otimes_{\O_E}k\ \ .
\end{eqnarray*}

Here we used (\ref{2.8}) in the first line and the previous proposition in
the last line. Now let ${\bf r}$ be arbitrary and let ${\bf t}={\bf
r}^{\vee}$ be the dual partition to ${\bf r}$, i.e.

\begin{eq}\label{3.13}
t_1=\# \{ \varphi; r_{\varphi}\geq 1\},\  t_2=\#\{ \varphi;r_{\varphi}\geq
2\}\ ,\ \hbox{etc.}
\end{eq}
By (\ref{2.2}) the partition ${\bf t}$ lies in ${\cal S}^0(r,e,d)$. Hence
${\bf t}\leq {\bf s}_{\rm max}(r,e)$ with equality only if ${\bf r}={\bf
r}_{\rm min}(r,e)$. Hence, if ${\bf r}\neq {\bf r}_{\rm min}(r,e)$, we
have
\begin{eq}\label{3.14}
{\rm dim}\ M_{\bf t} < {\rm dim}\  M_{{\bf s}_{\rm max}}={\rm dim}\
\overline M\ \ .
\end{eq}
Now
\begin{eqnarray*}
{\rm dim}\ M_{\bf t} &=& \langle {\bf t}, 2\varrho\rangle
=t_1(d-1)+t_2(d-3)+\ldots + t_d(1-d)\\ &=& d\cdot\sum^d_i t_i-\sum^d_{i=1}
(2i-1) t_i\ \ .
\end{eqnarray*}
But $\sum_i t_i =\sum_\varphi r_{\varphi} =r$. The second sum on the right
hand side can be written as a sum of contributions of each $\varphi$. Each
fixed $\varphi$ contributes $\sum^{r_{\varphi}}_{j=1}
(2j-1)=r^2_{\varphi}$. Hence

\begin{eq}\label{3.15}
{\rm dim}\ M_{\bf t} =dr-\sum_\varphi r^2_{\varphi} ={\rm dim}\
M(\Lambda,{\bf r})\otimes_{\O_E}E\ \ .
\end{eq}
Taking into account (\ref{3.14}), the result follows.
\end{Proof}\endproof

\begin{cor}\label{cor3.3}
If ${\bf r}\neq {\bf r}_{\rm min}(r,e)$, the corresponding standard model
$M$ is not flat over $\Spec\O_E$.
\endproof
\end{cor}

There remains the question whether if ${\bf r}={\bf r}_{\rm min}(r,e)$,
the corresponding standard model is flat over $\O_E$; we will return to it
in section \ref{flat}.

We end this section with the following remark. Among all strata of
$\overline M$ enumerated by ${\cal S}^0(r,e,d)$ there is a unique minimal
one,
\begin{eq}\label{3.16}
{\bf s}_{\rm min} ={\bf s}_{\rm min}(r,d)= ((u+1)^j, u^{d-j})\ \ .
\end{eq}
Here we have written $r=u\cdot d+j,\ 0\leq j<d$. The corresponding stratum
$M_{{\bf s}_{\rm min}}$ is closed and lies in the closure of any other
stratum.

\medskip

\section{Relation to the nilpotent variety} \label{nilpotent}

\setcounter{equation}{0}

In this section we describe the connection of the standard models for
$GL_d$ with the nilpotent variety. Consider the $GL_r$-torsor
$
\pi: \widetilde M\lo M
$
where
\begin{eq}\label{4.1}
\widetilde M(S)=\{ (\F, \psi)\ ;\ \F\in  M(S), \psi: \F\buildrel\simeq
\over\longrightarrow \O_S^r\}
\end{eq}
and $GL_r$ acts via $(\F, \psi)\mapsto (\F, \gamma\cdot\psi)$. Then $\Gg$
acts on $\widetilde M$ via $$(g, ({\cal F}, \psi))\longmapsto (g({\cal
F}), \psi\cdot g^{-1})\ \ ,$$ and the morphism $\pi$ is equivariant.

Note that $\prod_{\varphi}(T-\varphi ( \pi))^{r_{\varphi}}$ has
coefficients in $\O_E$; let us define a scheme $N=N({\bf r})$ over
$\Spec\O_E$ via the functor which to the $\O_E$-algebra $R$ associates the
set

\begin{eq}\label{4.2} \{A\ {\in}\ M_{r{\times} r}(R) ;\
{\rm det}(T\cdot I-A)\ {\equiv}\
\Pi_{\varphi}(T-\varphi(\pi))^{r_{\varphi}}, \ Q(A){=}0\}.
\end{eq}

The group scheme $GL_r$ acts on $N$ via conjugation
$A\mapsto \gamma\cdot A\cdot \gamma^{-1}$. There is a
morphism
\begin{eq}\label{4.3}
\phi: \widetilde M\ \lo\ N,\ \ \ \phi((\F, \psi))= \psi(\pi|\F)\psi^{-1}\
\ ,
\end{eq}
which is $GL_r$-equivariant. We have for $g\in {\cal G}$,
$$
\phi(g\cdot (\F, \psi))=\psi g^{-1}\cdot (\pi|g(\F))\cdot
(\psi g^{-1})^{-1}=
\psi\cdot (\pi|\F)\cdot \psi^{-1}=\phi((\F, \psi))
$$
since $g$ commutes with $\pi$. Therefore, the morphism $\phi$ is
$\Gg$-equivariant with trivial $\Gg$-action on the target $N$. We
therefore obtain a diagram of morphisms

\begin{eq}\label{4.4}
M\buildrel\pi\over\longleftarrow \widetilde M
\buildrel\phi\over\longrightarrow N\ \ ,
\end{eq}
which is equivariant for the action of $\Gg\times GL_r$, where the
first factor acts trivially on the right hand target and the second
factor acts trivially on the left hand target.

Since $Q(T)=\prod_{\varphi}(T-\varphi(\pi))$ with all elements
$\varphi(\pi)$ pairwise distinct, we see that the generic fiber of $N$
consists of all {\it semisimple} matrices with eigenvalues $\varphi(\pi)$
with multiplicity $r_{\varphi}$. It follows that $GL_r$ acts transitively
on $N\otimes_{{\cal O}_E}E$, which is smooth of dimension
$r^2-\Sigma_{\varphi}r^2_{\varphi}$ over $\Spec E$. The special fiber
$\overline N$ of $N$ is the subscheme of ${\rm Mat}_{r\times r}$ defined
by the equations

\begin{eq}\label{4.5}
\overline N=\{ A\in {\rm{Mat}}_{r\times r};\ A^e=0,\ {\rm
det}(T\cdot I-A)\equiv T^r\}\ \ .
\end{eq}
 This is a closed subscheme of the variety of nilpotent $r\times r$
matrices (which is defined by ${\rm det}(T\cdot I-A)\equiv T^r$).

The following result exhibits a close connection between the standard
model and the nilpotent variety.

\medskip

\begin{thm} \label{phismooth}
The morphism $\phi: \widetilde M\ \lon\ N$ is smooth of relative dimension
rd.
\end{thm}

\medskip

\begin{Proof}  Let $\Mod=\Mod(\O_F, {\bf r})$ be the algebraic stack
over $\Spec\O_E$ given by the fibered category of
$\O_F\otimes_{\O_{F_0}}\O_S$-modules ${\cal{F}}$ which are locally free
$\O_S$-modules of rank $r$ (comp. [LMoB] 3.4.4, 4.6.2.1) and for which $$
{\rm det}(T\cdot I-\pi|{\cal{F}})\equiv\prod_\phi(T-\phi(\pi))^{r_\phi}\
$$ as polynomials.  There is an isomorphism $$ \Mod\ \simeq\ [N/GL_r]\ $$
(where the quotient stack is for the conjugation action) given by
$\Gg\mapsto$ the conjugation $GL_r$-torsor of matrices $A$ giving the
action of $\pi$ on $\Gg$. By the definition of the quotient stack the
diagram (\ref{4.4}) corresponds to a morphism
\begin{eq}\label{4.50a}
\bar\phi: M\longrightarrow [N/GL_r]\ \ .
\end{eq}
The morphism $\phi$ is smooth if and only if $\bar\phi$ is a smooth
morphism of algebraic stacks. Under the identification above the morphism
$\bar\phi$ becomes $$ \bar\phi: M\to \Mod\ \ ;\quad (\F\subset
\Lambda\otimes_{\O_{F_0}}\O_S)\mapsto \F\ . $$ Now consider the morphism
$\phi^*: \widetilde M\to N$ obtained by composing $\phi$ with the
automorphism of $N$ given by $A\mapsto {^tA}$; $\phi$ is smooth if and
only if $\phi^*$ is smooth. Set $\F^*:=Hom_{\O_S}(\F, \O_S)$ which is also
naturally an $\O_F\otimes_{\O_{F_0}}\O_S$-module. By regarding $\widetilde
M$ as the $GL_r$-torsor over $M$ giving the $\O_S^r$-trivializations of
the dual $\F^*$ we see that $\phi^*$ descends to $$ \bar\phi^*: M\to \Mod\
\ ;\quad (\F\subset \Lambda\otimes_{\O_{F_0}}\O_S)\mapsto \F^*\ . $$ As
before, it is enough to show that $\bar\phi^*$ is smooth.

The stack $\Mod$ supports the universal
$\O_F\otimes_{\O_{F_0}}{\O_{\Mod}}$-module ${\cal{F}}$. Let us consider
the $\O_{\Mod}$-module $$ {\cal
T}=Hom_{\O_F\otimes_{\O_{F_0}}\O_{\Mod}}(\Lambda\otimes_{\O_{F_0}}\O_{\Mod},
{\cal F}) \simeq {\cal F}^{\oplus d}\ . $$ This defines a vector bundle
${\mathfrak V}({\cal T}^*)$ over $\Mod$ (see [LMoB] 14.2.6). The structure
morphism ${\mathfrak V}({\cal T}^*)\to \Mod$ is representable and smooth.
Its effect on objects is given by $$ ({\cal H}, \
\Lambda\otimes_{\O_{F_0}}\O_S\buildrel f\over\to {\cal H})\mapsto {\cal
H}\ . $$

Consider now the perfect $\O_{F_0}$-bilinear pairing $$ (\ ,\
):\O_F\times\O_F\to \O_{F_0}\ \ ;\quad (x,y)={\rm
Tr}_{F/F_0}(\delta^{-1}xy) $$ where $\delta$ is an $\O_F$-generator of the
different ${\cal D}_{F/F_0}$. This pairing gives an $\O_F$-module
isomorphism $\O_F^*=Hom_{\O_{F_0}}(\O_F, \O_{F_0})\simeq \O_F$. Now choose
an $\O_F$-module isomorphism $\Lambda\simeq\O_F^d$. We obtain functorial
$\O_F\otimes_{\O_{F_0}}\O_S$-module isomorphisms
$(\Lambda\otimes_{\O_{F_0}}\O_S)^*\simeq \Lambda\otimes_{\O_{F_0}}\O_S$
and a morphism $i: M\to {\mathfrak V}({\cal T}^*)$ given by $$
({\F}\subset \Lambda\otimes_{\O_{F_0}}\O_S)\ \mapsto\ (\F^*,\
\Lambda\otimes_{\O_{F_0}}\O_S\simeq (\Lambda\otimes_{\O_{F_0}}\O_S)^*\to
\F^*). $$ We therefore see that $i$ is representable and an open immersion
(it gives an isomorphism between $M$ and the open substack of ${\mathfrak
V}({\cal T}^*)$ given by the full subcategory of objects for which the
morphism $f$ is surjective). The morphism $\bar\phi^*$ is the composition
$$ M\to {\mathfrak V}({\cal T}^*)\to \Mod$$ and therefore is smooth. The
relative dimension of $\phi$ is equal to the relative dimension of
$\bar\phi$; this, in turn, is equal to the relative dimension of the
composition above. However, this is the same as the relative dimension of
the vector bundle ${\mathfrak V}({\cal T}^*)\to \Mod$ which is equal to
$rd$.
\end{Proof}

\begin{Remarks}\label{remarks4.5}
{\rm (i) The proof of Theorem \ref{phismooth} given above follows ideas
which appear in the paper [FGKV] (see especially loc.cit. \S 4.2) and were
brought to the attention of the authors by D. Gaitsgory. A previous
version of the paper contained a more complicated proof which used
explicit matrix calculations to describe affine charts for the scheme $M$.
\smallskip

\par\noindent
 (ii) Let us consider the conjugation action of $GL_r$ on the
special fibre $\overline N$. The orbits of this action are
parametrized by
\begin{eq}\label{4.48}
{\cal S}(r,e)=\{ {\bf s}=(s_1\geq \ldots \geq s_r);\ e\geq s_1,\ s_r\geq
0,\ \Sigma_i s_i=r\}\ \ .
\end{eq}
We denote the corresponding orbit by $N_{\bf s}$. Again $N_{{\bf s}'}$
lies in the closure of $N_{\bf s}$ if and only if ${\bf s}'\leq {\bf s}$.
Obviously ${\cal S}^0(r,e,d)\subset{\cal S}(r,e)$ and the ${\cal G}\times
GL_r$-equivariance of the diagram (\ref{4.4}) shows that for ${\bf s}\in
{\cal S}^0(r,e,d)$,
\begin{eq}\label{4.49}
\phi(\pi^{-1}(M_{\bf s}))=N_{\bf s}\ \ .
\end{eq}
In particular, the image of $\tilde M\otimes_{\O_E}k$ under $\phi$ is the
union of orbits corresponding to ${\bf s}\in {\cal S}^0(r,e,d)$ and the
diagram (\ref{4.4}) induces an injection of the set of ${\cal G}$-orbits
in $\overline{M}$ into the set of $GL_r$-orbits in $\overline N$. It is
easy to see that the complement of ${\cal S}^0(r,e,d)$ in ${\cal S}(r,e)$
is closed under the partial order $\leq$ on ${\cal S}(r,e)$, i.e.\
corresponds to a closed subset of $\overline N$.
\smallskip

\par\noindent
(iii) The dimension of $N_{\bf s}$ is given by the formula
\begin{eq}\label{4.50}
{\rm dim}\ N_{\bf s} =r^2-\sum_{i=1}^e r_i^2\ \ ,
\end{eq}
where $r_1\geq \ldots \geq r_e\geq 0$ is the dual partition to ${\bf s}$.
\par\noindent
For ${\bf s}$ in ${\cal S}^0(r,e,d)$, this formula is compatible with the
one of (\ref{3.14}) via Theorem \ref{phismooth}: we have $$\langle {\bf
s}, 2\varrho\rangle +r^2-rd=r^2-\sum_{i=1}^e r_i^2\ \ ,$$ comp.\
(\ref{3.15}) (which is ``dual'' to (\ref{4.50})).
\smallskip

(iv) The equivariance of (\ref{4.4}) can be rephrased (via descent theory for
${\cal G}$-torsors) by saying that the morphism (\ref{4.50a}) factors
through a morphism of algebraic stacks
\begin{eq}\label{4.50b}
[M/{\cal G}]\longrightarrow [N/GL_r]\ \ .
\end{eq}
Note that both these stacks have only finitely many points. (The set of
points in the special fiber of $[M/{\cal G}]$ resp.\ $[N/GL_r]$ is ${\cal
S}^0(r,e,d)$ resp.\ ${\cal S}(r,e)$.)}
\end{Remarks}

\begin{cor}\label{cor4.6}
Let $r_{\varphi}\leq 1,\ \forall\varphi$. Then the standard model for
$GL_d$ corresponding to ${\bf r}$ is flat over $\Spec\O_E$, with special
fiber a normal complete intersection variety.
\end{cor}

\begin{Proof}
In this case the first condition $A^e=0$ in the definition (\ref{4.5}) of
$\overline N$ is a consequence of the second condition. Hence $\overline
N$ is the variety of nilpotent matrices, which is a reduced and
irreducible, normal and complete intersection variety. On the other hand,
obviously ${\bf r}={\bf r}_{\rm min}$, hence by Proposition \ref{prop3.2}
we have ${\rm dim}\ N\otimes_{\O_E}E={\rm dim}\ N\otimes_{\O_E}k$. hence
the generic point of $\overline N$ is the specialization of a point of
$N\otimes_{\O_E}E$. Since $N\otimes_{\O_E}k$ is reduced, the flatness of
$N$ follows from EGA IV 3.4.6.1. By Theorem \ref{phismooth} this implies
the corresponding assertions for $M$.\endproof
\end{Proof}

\begin{Remark}\label{remark4.7} {\rm Let us fix an isomorphism $\Lambda\simeq\O_F^d$ and
suppose that $r=d=e$. Write $$ \Lambda\otimes_{\O_{F_0}}\O_S\simeq
(\O_F\otimes_{\O_{F_0}}\O_S)^r= 1\cdot\O_S^r\oplus\pi\cdot
\O_S^r\oplus\cdots \oplus\pi^{r-1}\cdot\O_S^r\ . $$ Now let us consider
$\O_F\otimes_{\O_{F_0}}\O_S$-submodules
$\F\subset\Lambda\otimes_{\O_{F_0}}\O_S $ which are locally free
$\O_S$-locally direct summands and are such that the composition $$ F\ :\
\F\subset \Lambda\otimes_{\O_{F_0}}\O_S \buildrel pr\over\longrightarrow
\pi^{r-1}\cdot\O_S^d $$ is a surjection (and therefore an isomorphism).
The inverse of the isomorphism $F$ can be written $$ F^{-1}(\pi^{r-1}\cdot
v)=(1\cdot f_{r-1}(v),\ \pi\cdot f_{r-2}(v),\ \ldots ,\ \pi^{r-1}\cdot
f_0(v)) $$ where $f_i:\O_S^r\to\O_S^r$, $i=0,\cdots, r-1$, are
$\O_S$-linear homomorphisms and $f_0={\rm Id}$. Recall $Q(T)=
T^e+\sum_{k=0}^{e-1}b_kT^k$. The condition that $\F$ is stable under
multiplication by $\pi$ translates to $$ f_{k+1}-b_{r-k-1}=f_{k}\cdot
(f_1-b_{r-1}\cdot f_0), \ \ k=1,\ldots, r-2\ , $$ $$
 -b_0=f_{r-1}\cdot (f_1-b_{r-1}\cdot f_0)=0\ .
$$ Therefore, all the $f_k$ are determined by $f_1$ and we have
$Q(f_1-b_{r-1}\cdot {\rm Id})=0$. We also have $F\cdot (\pi|\F)\cdot
F^{-1}=f_1-b_{r-1}\cdot {\rm Id}$. Using these facts, we see that after
choosing a basis of $\O_S^r$, the modules $\F$ for which the composition
$F$ is an isomorphism are in 1--1 correspondence with $r\times r$-matrices
$A$ which satisfy $Q(A)=0$. Let $U$ denote the open subscheme of $M$ whose
$S$-points correspond to modules $\F$ for which the homomorphism $F$ above
is an isomorphism. In terms of the diagram (\ref{4.4}), the above implies
that, in this case, there is a section $s: N\to \widetilde M$ to the
morphism $\phi:\widetilde M\to N$ such that the composition $\pi\circ s:
N\to M$ is an open immersion identifying $N$ with $U$. The smoothness of
the morphism $\phi$ then amounts to the smoothness of the conjugation
action morphism $N\times GL_r\to N$. Also, since $M$ is projective, this
shows that, in this case, the local model $M$ may be considered as a
relative compactification over ${\rm Spec}\ {\cal O}_E$of the variety $N$.
This last result for the special fibers is precisely the scheme-theoretic
version of Lusztig's result [L], section 2. Hence the special fiber of any
standard model for $GL_d$, for which $r=d$ and $e=d$ (they are all
identical), may be considered as a compactification of the nilpotent
variety. We return in section \ref{section5A} to the consequences of
Theorem \ref{phismooth} for the special fibers.}
\end{Remark}

\medskip

\section{The canonical flat model} \label{flat}

\setcounter{equation}{0}

We have seen in Corollary \ref{cor3.3} that a standard model is rarely
flat over ${\Spec}\O_E$. By Theorem \ref{phismooth} the same can be said
of the scheme $N$. In this section we will first show that the flat scheme
theoretic closure of the generic fiber $N\otimes_{\O_E}E$ in $N$ has good
singularities. The idea is to use a variant of the Springer resolution of
the nilpotent variety, as also in the work of Eisenbud and Saltman
([E-S]).

Recall our notations from the beginning of section \ref{GLd}. Let $K$ be
the Galois hull of $F$ inside $F^{\rm sep}_0$. Let us order the different
embeddings $\phi:F\rightarrow K$. Then we can write
\begin{eq}\label{5.1}
P(T)=\prod_{\phi}(T-\phi(\pi))^{r_{\phi}}=
\prod_{i=1}^e(T-a_i)^{r_i}\ ,\ Q(T)=\prod_{i=1}^e(T-a_i)\ .
\end{eq}
Here $\phi(\pi)=a_i\in \O_K$ are distinct roots.

Let us set $ n_k=\sum_{i=1}^kr_i $,  for $1\leq k\leq e$. Let {\boldmath
$\cal F$} be the scheme which classifies flags:

\begin{eq}\label{5.2}
(0)={\cal F}_{e}\subset {\cal F}_{e-1}\subset \cdots \subset {\cal
F}_{0}=\O_S^r
\end{eq}
where ${\cal F}_k$ is locally on $S$ a direct summand of $\O_S^r$ of
corank $n_k$. Following [E-S], we consider the subscheme
$\N$ of $(({\rm Mat}_{r\times r})\times \hbox{{\boldmath $\cal
F$}})_{\O_K}$ classifying pairs $(A, \{{\cal F}_{\bullet}\})$ such
that
\begin{eq}\label{5.3}
(A-a_k\cdot I)\cdot {\cal F}_{k-1}\subset {\cal F}_{k},\quad 1\leq
k\leq e\ .
\end{eq}
The scheme $\N$ supports an action $GL_r$ by $g\cdot (A, \{{\cal
F}_i\})=(gAg^{-1}, \{g({\cal F}_i)\})$.

Obviously this is a variant of the Grothendieck-Springer
construction. It differs from the original in two aspects: we
consider partial flags instead of complete flags, and we fix the
(generalized) eigenvalues $a_i$ of $A$.

\begin{lemma} \label{natu}

i) $\N$ is smooth over $\Spec\O_K$.

ii) There is a projective $GL_r$-equivariant morphism $\mu:\N\ \lo\
N\otimes_{\O_E} \O_K$, given by $(A, \{{\cal F}_i\})\mapsto A$.

iii) The morphism $\mu_{{\Spec}K}$ is an isomorphism between the generic
fibers $\N\otimes_{\O_K}K$ and $N\otimes_{\O_E}K$.
\end{lemma}

\begin{Proof} Part (i) follows from the fact that the projection to the
second factor is smooth

\begin{eq}\label{5.4}
\N\longrightarrow \hbox{{\boldmath $\cal F$}}\ ,
\end{eq}
comp.\ [E-S], p.\ 190 (the fiber over $\{{\cal F}_i\}$ can be identified
with the cotangent space of $\bcf$ at $\{{\cal F}_i\}$). Now suppose that
$A$ is in $M_{r\times r}(R)$ for some $\O_K$-algebra $R$ and that locally
on Spec $R$ there is a filtration $\{{\cal F}_{\bullet}\}$ of $R^r$ as
described above. Then the characteristic polynomial of $A$ is equal to
$P(T)$ and we have

\begin{eq}\label{5.5}
\prod_{k=1}^N(A-a_k\cdot I)=0\ \in\ M_{r\times r}(R).
\end{eq}
This implies $Q(A)=0\in M_{r\times r}(R)$, cf.\ (\ref{5.1}). This shows
that the natural morphism $\N\ \lo\ {({\rm Mat}_{r\times r})}_{\O_K}$
factors through $N\otimes_{\O_E}\O_K$; the claim (ii) follows. Now
$N\otimes_{\O_E} K$ consists of all semisimple matrices with eigenvalues
$a_i$ with multiplicity $r_i$ (comp.\ remarks before (\ref{4.5})). Hence
the fiber of $\mu_{{\rm Spec}\ K}$ over $A$ is the filtration ${\cal
F}_{\bullet}$ associated to the eigenspace decomposition corresponding to
$A$, i.e., is uniquely determined by $A$.
\endproof
\end{Proof}

\medskip

Now set $N'={\rm Spec}({\mu}_*(\O_{\N}))$ and denote by $\mu(\N)$ the
scheme-theoretic image of $\mu:\N\lo N\otimes_{\O_E}\O_K$; this is a
reduced closed subscheme of $N\otimes_{\O_E}\O_K$ since $\N$ is smooth
(therefore reduced) and $\mu$ is proper. The scheme $N'$ supports an
action of $GL_r$; there is a natural $GL_r$-equivariant morphism $q':
N'\ \lo\ N\otimes_{\O_E}\O_K$. The morphism $q'$ is finite and factors as
follows

\begin{eq}\label{5.6}
q': N'\lo \ \mu(\N)\lo \ N\otimes_{\O_E}\O_K.
\end{eq}
Lemma \ref{natu} (iii) implies that $q'_{\Spec K}$ is an isomorphism.

Let $k'$ be the residue field of $\O_K$.

\begin{prop} \label{normal}

a) The special fiber $\overline N'=N'\otimes_{\O_K}k'$ is normal
and has rational singularities.

b) The scheme $N'$ is normal, Cohen-Macaulay and flat over $\Spec\O_K$.

c) $N'=\mu(\N)$.

d) $N'$ is the scheme theoretic closure of $N\otimes_{\O_E}K$ in
$N\otimes_{\O_E}\O_K$. Its special fiber is the reduced closure of the
orbit $N_{\bf{t}}$, where ${\bf t} ={\bf r}^{\vee}$ is the dual
partition to ${\bf r}=(r_i)$.
\end{prop}

\begin{Proof} This follows closely
the arguments of [E-S] (see p. 190-192, proof of Theorem 2.1) with new
input the results of Mehta-van der Kallen ([M-vdK]). They show (using
Frobenius splitting) that the closure of the orbit of a nilpotent matrix
is normal and Cohen-Macaulay also in positive characteristic. For the
duration of this proof, we will denote by $\varpi$ a uniformizer of
$\O_K$. We will first consider the situation over $\Spec k'$ and use a bar
to denote base change from $\Spec\O_K$ to $\Spec k'$. Consider the
morphism $\overline\mu: \ \overline \N\ \lo\ \overline N$ and the scheme
$\Spec(\overline\mu_*(\O_{\overline\N}))\ \lo\ \overline N$. The morphism
$\overline\mu$ factors as
\begin{eq}\label{5.7}
\overline\mu: \ \overline \N\ \lo\
\Spec(\overline\mu_*(\O_{\overline\N}))\ \lo\ \overline\mu(\overline\N)\
\lo\ \overline N\ \ ,
\end{eq}
where again $\overline\mu(\overline\N)$ denotes the scheme theoretic image
of $\overline\mu$. The morphism $\overline\N\ \lo\
\overline\mu(\overline\N)$ is one to one on the open subset of
$\overline\N$ of those $(A, \{\F_{\bullet}\})$ such that $A{\cal
F}_{k-1}={\cal F}_k$, $k=1,\ldots, e$ and so it is birational; therefore
the morphism $\Spec(\overline\mu_*(\O_{\overline \N}))\ \lo\
\overline\mu(\overline\N)$ is finite and birational. Since $\overline\N$
is reduced, $\overline\mu(\overline\N)$ is also reduced. As in [E-S], we
see that $\overline\mu(\overline\N)\subset {({\rm Mat}_{r\times r})}_{k'}$
is the reduced closure of the conjugation orbit $N_{\bf t}$ of the Jordan
form for the dual partition ${\bf t}={\bf r}^{\vee}$. By [M-vdK], the
closure of $N_{\bf t}$ is normal and has rational singularities. We
conclude that $\Spec(\overline\mu_*(\O_{\overline\N}))\ \lo\
\overline\mu(\overline\N)$ is an isomorphism which we use to identify
these two schemes. We will first show that
$\Spec(\overline\mu_*(\O_{\overline\N}))=\overline\mu(\overline\N)$
actually gives the special fiber of $N'$. The statement (a) then follows
from the results of Mehta-van der Kallen. The cohomology exact sequence
obtained by applying $\mu_*$ to $$ 0\ \lo\ \O_{\N}\ \buildrel \varpi \over
\lo\ \O_{\N}\ \lo\ \O_{\overline\N}\ \lo\ 0 $$ gives an injective
homomorphism $$ \O_{N'}/\varpi\O_{N'}=\mu_*(\O_\N)/\varpi\mu_*(\O_\N)\
\lo\ \overline\mu_*(\O_{\overline\N}) $$ and it is enough to show that
this is an isomorphism. There is a commutative diagram $$
\matrix{\O_{\mu(\N)}/\varpi \O_{\mu(\N)}&\lo
&\O_{\overline\mu(\overline\N)}\cr \downarrow&&\downarrow\cr
\O_{N'}/\varpi\O_{N'} &\lo &\overline\mu_*(\O_{\overline\N}).\cr} $$ From
the definition of the scheme theoretic image, the upper horizontal
homomorphism is an isomorphism. On the other hand, we have seen above that
the right vertical homomorphism is an isomorphism, whence the claim. We
will now show that $N'\ \lo\ \mu(\N)$ is an isomorphism; this will
establish (c). From the above it follows that the homomorphism $$
\O_{\mu(\N)}/\varpi \O_{\mu(\N)}\ \lo\ \O_{N'}/\varpi\O_{N'} $$ is
surjective. Since $\mu$ is proper, $\O_{N'}=\mu_*(\O_{\N})$ is finite over
$\O_{\mu(\N)}$. Therefore, using Nakayama's lemma, we conclude that
$\O_{\mu(\N)}\ \lo\ \O_{N'}$ is surjective locally over all points at the
special fiber. Since $N'\otimes_{\O_K}K\ \lo\ \mu(\N)\otimes_{\O_K}K$ is
an isomorphism it follows that $\O_{\mu(\N)}\ \lo\ \O_{N'}$ is surjective;
it now follows from the definition of the scheme-theoretic image that
$\O_{\mu(\N)}\ \lo\ \O_{N'}$ is an isomorphism. This shows (c).

Now let us show part (b). By (a) the special fiber $\overline N'$ is
normal and Cohen-Macaulay; in fact, it has dimension $r^2-\sum_i r^2_i$.
This is equal to the dimension of the generic fiber
$N\otimes_{\O_E}K=N'\otimes_{\O_K}K$. As a result the special fiber is
reduced and its unique generic point lifts to the generic fiber; this
implies that $N'$ is flat over $\Spec\O_K$ (EGA IV 3.4.6.1). Since
$\overline N'$ is Cohen-Macaulay and $N'\ \lo\ \Spec\O_K$ is flat, $N'$ is
Cohen-Macaulay. Now $N'\otimes_{\O_K}K=N\otimes_{\O_E}K$ is smooth and
$\overline N'$ generically smooth; this shows that $N'$ is regular in
codimensions $0$ and $1$ and therefore, by Serre's criterion, normal.

Finally, since $N'=\mu(\N)\subset N\otimes_{\O_E}\O_K$ with identical
generic fibers, and since $N'$ is flat over $\Spec\O_K$, $N'$ is the
(flat) scheme theoretic closure of $N\otimes_{\O_E}K$ in $N\otimes_{\O_E}\O_K$. This shows
(d).\endproof
\end{Proof}

\medskip

\begin{prop}\label{prop5.4}
Let $N^{\rm loc}$ be the (flat) scheme theoretic closure of
$N\otimes_{\O_E}E$ in $N$. Then the scheme $N^{\rm loc}$ is normal and
Cohen-Macaulay. Its special fiber is the reduced closure of the orbit
$N_{\bf t}$ with ${\bf t}={\bf r}^{\vee}$ the dual partition to ${\bf
r}=(r_{\varphi})$. The special fiber is normal with rational
singularities.
\end{prop}

\begin{Proof}
Denote by $N''\subset N$ the scheme theoretic image of the finite
composite morphism $N'\ \lo\ N\otimes_{\O_E}\O_K\ \lo\ N$. This is a
$GL_r$-equivariant closed subscheme of $N$. We have $\O_{N''}\subset
\O_{N'}$ and so since $N'$ is flat over $\Spec\O_K$, we conclude that
$N''$ is also flat over $\Spec\O_E$. We have
$N''\otimes_{\O_E}E=N\otimes_{\O_E}E$; hence $N''$ is the flat scheme
theoretic closure of $N\otimes_{\O_E}E$ in $N$, that is $N''=N^{\rm loc}$.
The base change $N^{\rm loc}\otimes_{\O_E}\O_K$ is a closed subscheme of
$N\otimes_{\O_E}\O_K$ which is flat over $\Spec\O_K$. Since we have
$N^{\rm loc}\otimes_{\O_E}K=N\otimes_{\O_E}K$ by Proposition \ref{normal}
(d) we have
\begin{eq}\label{5.8}
N^{\rm loc}\otimes_{\O_E}\O_K= N'\subset N\otimes_{\O_E}\O_K\ \ .
\end{eq}
By Proposition \ref{normal} (b) $N'$ is normal. Therefore, $N^{\rm
loc}=N''/{\rm Gal}(K/E)$ is also normal. Denoting by $\overline N^{\rm
loc}$ the special fiber of $N^{\rm loc}$, we have $$ \overline N^{\rm
loc}\otimes_kk'=\overline N'\  . $$ Hence the remaining assertions also
follow from Proposition \ref{normal}.\endproof
\end{Proof}

\medskip
We note that we can define analogues of ${\cal N}$ for $M$ and $\widetilde
M$. Namely, we consider the $\O_K$-scheme ${\cal M}$, which for an
$\O_K$-scheme $S$ classifies the filtrations of $\O_S$-submodules

\begin{eq}\label{5.8a}
\{ (0)= {\cal F}_e\subset {\cal F}_{e-1}\subset\ldots\subset {\cal F}_0=
{\cal F}\subset \Lambda\otimes_{\O_{F_0}}\O_S\}\ ,
\end{eq}
where ${\cal F}\in M^{\rm loc}(S)$ and where ${\cal F}_k$ is locally on
$S$ a direct summand of $\cal F$ of corank $n_k=\sum_{i=1}^kr_i$ such that

\begin{eq}\label{5.8b}
(\pi-a_k\cdot {\rm Id})\cdot {\cal F}_{k-1}\subset {\cal F}_k\  ,\ \
k=1,\ldots, e\ .
\end{eq}

Similarly we define $\widetilde{\cal M}$ by fixing in addition an
isomorphism $\psi:{\cal F}\to \O_S^r$. We thus obtain a diagram with
cartesian squares in which the vertical morphisms are projective with
source a smooth $\O_K$-scheme,

\begin{eq}\label{5.8c}
\matrix{ {\cal M} & \longleftarrow & \widetilde{\cal M} & \longrightarrow
& {\cal N} \cr \big\downarrow && \big\downarrow && \big\downarrow \cr M
\otimes_{\O_E}\O_K &\longleftarrow & \widetilde M \otimes_{\O_E}\O_K &
\longrightarrow & N\otimes_{\O_E}\O_K & . \cr}
\end{eq}

We now use Theorem \ref{phismooth} to transfer the previous results to the
standard models.

\begin{thm}  \label{mcan}
Let $M^{\rm loc}$ be the (flat) scheme theoretic closure of $M
\otimes_{\O_E}E$ in $M $. Then

(i) $M^{\rm loc}$ is normal and Cohen-Macaulay.

(ii) The special fiber $\overline M^{\rm loc}$ is reduced, normal with
rational singularities, and we have $$ \overline M^{\rm
loc}=\bigsqcup_{{\bf s}\leq {\bf t}} M^{\rm loc}_{\bf s}. $$

(iii) There is a diagram with cartesian squares of ${\cal G}\times
GL_r$-equivariant morphisms in which the horizontal morphisms are smooth
and the vertical morphisms are closed embeddings, $$\matrix{ M^{\rm loc} &
\buildrel \pi^{\rm loc}\over\longleftarrow & \widetilde M^{\rm loc} &
\buildrel \phi^{\rm loc}\over\longrightarrow & N^{\rm loc} \cr
\big\downarrow && \big\downarrow && \big\downarrow \cr M &
\buildrel\pi\over\longleftarrow & \widetilde M &
\buildrel\phi\over\longrightarrow & N & .&&
\endproof
\cr}$$

\end{thm}

The $\O_E$-scheme $M^{\rm loc}$ is called the {\it canonical model}
associated to the standard model for $GL_d$ corresponding to ${\bf
r}=(r_i)$.

\begin{Remark} \label{mathieu}
{\rm The use of Theorem \ref{phismooth} in proving parts (i) and (ii) of
the above theorem is to enable us to appeal to the results of Mehta-van
der Kallen on nilpotent orbit closures. An alternative approach which
would not appeal to Theorem \ref{phismooth} might be obtained by studying
directly the composite morphism $$ {\cal M}\ \lo\ M\otimes_{\O_E}\O_K\
\lo\ M $$ and using the theory of generalized Schubert varieties in the
affine Grassmannian for $GL_d$. One can show directly that ${\cal M}$ is
smooth over $\Spec\O_K$ (comp. the proof of Lemma \ref{natu}); the main
point is that the special fiber $\overline{\cal M}$ can be written as a
composite of smooth fibrations with fibers Grassmannian varieties (indeed,
we may think of $\overline{\cal M}$ as a generalized
Demazure-Bott-Samelson variety). Then, the arguments in the proofs of
Propositions \ref{normal} and \ref{prop5.4} can be repeated to obtain a
direct proof of parts (i) and (ii) of Theorem \ref{mcan}, provided we know
that the reduced closure of the stratum $M_{\bf t}$ in the special fiber
$\overline{M}$ is normal and has rational singularities. By the discussion
in \S 3, this reduced closure is isomorphic to the reduced closure $X_{\bf
t}:=\overline\O_{\bf t}$ of the Schubert cell $\O_{\bf t}$ in the affine
Grassmannian (defined as in [BL], comp. \S3). It remains to show that the
generalized Schubert variety $X_{\bf t}$ is normal and has rational
singularities. Results of this type have been shown (in positive
characteristic) by Mathieu [Mat]. Unfortunately, it is not clear that in
positive characteristic the Schubert varieties he considers have the same
scheme structure as the $X_{\bf t}$ (see the remarks on p. 410 of [BL]).
Since this point is not cleared up, we use Theorem \ref{phismooth} to
deduce results on the affine Grassmannian from results on nilpotent orbit
closures, comp.\ section \ref{section5A} below. }
\end{Remark}

\smallskip

Suppose now that $(|\Gamma|, {\rm char}\ k)=1$. Under this hypothesis we
will show that a conjecture of de Concini and Procesi implies a rather
explicit description of the scheme $M^{\rm loc}$.

We follow [deC-P] \S  1. Let $x_1$, $\ldots$, $x_r$ be a set of
variables; for every pair of integers $t$, $h$ with $h\geq 0$,
$1\leq t\leq r$, we can consider the total symmetric function of
degree $h$ in the first $t$ variables; this is defined to be the
sum of all monomials in $x_1$, $\ldots$, $x_t$ of degree $h$ and
will be denoted by $S^t_h(x_i)$.  We also indicate by the symbols
$\sigma_h$ the elementary symmetric function of degree $h$ in the
variables $x_1$, $\ldots$, $x_r$ with the convention that
$\sigma_h=0$ if $h>r$. Write
\begin{eq}\label{5.9}
S_h^t(x_i)=\sum a_{(h_1,\cdots, h_t)}x_1^{h_1}\cdots x_{t}^{h_t}
\end{eq}
For $A\in {\rm Mat}_{r\times r}(R)={\rm End}(R^r)$, we set
\begin{eq}\label{5.10}
S^t_h(A)(e_{i_1}\wedge  \cdots \wedge e_{i_t})=
\sum a_{(h_1,\cdots, h_t)}A^{h_1}e_{i_1}\wedge \cdots \wedge A^{h_t}e_{i_t}
\end{eq}
Since $S^t_h(x_i)$ is symmetric this defines a $R$-linear operator
\begin{eq}\label{5.11}
S^t_h(A): \wedge^t(R^r)\ \lo\ \wedge^t(R^r)
\end{eq}
Now let us indicate by
\begin{eq}\label{5.12}
T^r-\sigma_1(A)T^{r-1}+\sigma_2(A)T^{r-2}-\cdots +(-1)^r\sigma_r(A)
\end{eq}
the characteristic polynomial ${\rm det}(T\cdot I-A)$ of $A$. For $t$, $h$
as above, we now define the following element of ${\rm End}(\wedge
^t(R^r))$: $$ F^t_h(A):=S^t_h(A)-\sigma_1(A)S^t_{h-1}(A)+\cdots
+(-1)^h\sigma_h(A) $$

For each function $\phil: \{1,2,\ldots, e\}\ \lo\ {\bf N}$ with
 $0\leq \phil(i)\leq r_i$, for all $1\leq i\leq e$, consider
\begin{eq}\label{5.13}
Q_{\phil}(t)=\prod^e_{i=1}(T-a_i)^{\phil(i)}
\end{eq}
(a divisor of the polynomial $P(T)$).

Let us consider the subscheme $N_0$ of ${\rm Mat}_{r\times r}$ over $\Spec\O_E$
which is defined by the equations given by the conditions
\begin{eq}\label{5.14}
{\rm det}(T\cdot I-A)\equiv P(T), {\rm \ and\ \ }
\
\sum_{\sigma\in \Gamma} F^t_h(A)\cdot \wedge^t(\sigma Q_\phil(A))=0,
\end{eq}
$$
{\rm for\ all\   }\phil,{\rm{\ \ and\ for}\ \ }t+h=r-\sum_{i,\phil(i)\neq 0}r_i+1,\ \ t\geq 1,\ h\geq 0.
$$
(it is obvious that the generators of the ideal defining $N_0$ have
all coefficients in $\O_E$). This is actually a closed subscheme of
$N$; indeed, consider the second set of equations for $\phil\equiv
1$, $t=1$, $h=0$. Since $F_0^t(A)=I$, we obtain
\begin{eq}\label{5.15}
|\Gamma|\cdot \prod_{i=1}^e (A-a_i\cdot I)=0\ .
\end{eq}
Since $(\vert\Gamma\vert, {\rm char}\ k)=1$, this equation implies
$Q(A)=0$. It is straightforward to see that $N_0$ is a
$GL_r$-invariant subscheme of $N$.

\begin{prop} \label{generic}
i) The schemes $N$ and $N_0$ over $\Spec\O_E$ have the same generic
fiber.

ii) The reduced special fiber of $N_0$ is equal to the reduced closure
of the conjugation orbit $N_{\bf t}$ of the Jordan form with
partition ${\bf t}={\bf r}^{\vee}$.
\end{prop}

\begin{Proof} We first show (i). By descent, we can check this after
base changing to $K$. Consider the diagonal $r\times r$
matrix $A_0$ in which the element $a_i$ appears with multiplicity $r_i$.
We have

\begin{eq}\label{5.16}
{\rm rank}(Q_{\phil}(A_0))=r-\sum_{i,\phil(i)\neq 0}r_i.
\end{eq}
It now follows from [deC-P], Proposition on p. 206,
that
\begin{eq}\label{5.17}
F^t_h(A_0)\cdot \wedge^t (Q_{\phil}(A_0))=0
\end{eq}
if $t+h\geq r-\sum_{i,\phil(i)\neq 0}r_i+1$ (comp.\ the proof of the
Theorem on p. 207 loc. cit). This shows that $A_0$ satisfies the
conditions defining ${N_0}\otimes_{\O_E}K$. Since $N\otimes_{\O_E}K$ is
the reduced $GL_r$-orbit of such a matrix and $N_0$ is
$GL_r$-equivariant, the result follows.

Now we prove (ii). Reducing the equations defining $N_0$ modulo the
maximal ideal of $\O_E$ gives the following equations:
\begin{eq}\label{5.18}
{\rm det}(T\cdot I-A)\equiv T^r, {\rm \ and\ \ } |\Gamma|\cdot
F^t_h(A)\cdot \wedge^t(A^{{\sum_i\phil(i)}})=0
\end{eq}
$$
{\rm for\ all\   }\phil,{\rm{\ \ and\ for}\ \ }
t+h=r-\sum_{i,\phil(i)\neq 0}r_i+1,\ \ t\geq 1,\ h\geq 0.
$$
By the definition of $F^t_h(A)$, the difference $F^t_h(A)-S^t_h(A)$
is in the ideal given by ${\rm det}(T\cdot I-A)\equiv T^r$. Using
this and taking $J={\rm supp}(f)\subset I=\{ 1,\ldots, e\}$ we see
that the above set of equations generates the same ideal as the
following one:
\begin{eq}\label{5.19}
{\rm det}(T\cdot I-A)\equiv T^r, {\rm \ and\ \ }
 S^t_h(A)\cdot \wedge^t(A^{|J|})=0
\end{eq}
$$ {\hbox{\rm for all subsets}\ }J\subset I, \ \ t+h=r-\sum_{i\in
J}r_i+1,\ \ t\geq 1,\ h\geq 0. $$ If the $r_i$ are arranged in decreasing
order, it is the same to consider
\begin{eq}\label{5.20}
{\rm det}(T\cdot I-A)\equiv T^r, {\rm \ and\ \ }
 S^t_h(A)\cdot \wedge^t(A^{k})=0
\end{eq}
$$ k=0,\ldots, e, \ \ t+h=r-n_k+1,\ \ t\geq 1,\ h\geq 0. $$ (recall our
notation $n_k=\sum_{i=1}^kr_i$, $n_0=0$). In fact, we could also omit the
first equation and write this as
\begin{eq}\label{5.21}
S^t_h(A)\cdot \wedge^t(A^{k})=0
\end{eq}
$$ \hbox{\ for\ } k=0,\ldots, e, \ \ t+h=r-n_k+1,\ \ t\geq 1,\ h\geq 0. $$
Indeed, by [deC-P], Lemma p. 206, the ideal generated by $S^{t}_h(A)$ for
$t+h=r+1$ is the same as the one generated by the coefficients of the
characteristic polynomial $\sigma_1(A)$, $\ldots$, $\sigma_r(A)$.

Let us denote by $\I$ the ideal of $k[a_{ij}]_{1\leq i,j\leq r}$ generated
by the equations (\ref{5.21}). Then it follows as in \hbox{[deC-P]},
Theorem on p. 207, that its radical ${\rm rad}(\I)$ defines the (reduced)
closure of the orbit $N_{\bf t}$ of the Jordan form for the dual
partition ${\bf t}={\bf r}^{\vee}$. This shows (ii).\endproof
\end{Proof}

\medskip

De Concini and Procesi conjecture (loc. cit.) that when
char $k=0$, we have ${\rm rad}(\I)=\I$. Assume that this is true
even if char $k>0$. Then

\begin{eq}\label{5.22}
N_0\otimes_{\O_E}k=\overline N_{\bf t}\ \ ,
\end{eq}
(i.e, the special fiber of $N_0$ is already reduced). Since by [M-vdK] the
reduced closure of $N_{\bf t}$ is normal, Cohen-Macaulay and has dimension
$r^2-\sum_ir_i^2$, an argument as in the proof of Proposition \ref{normal}
(ii) shows that $N_0$ is normal and flat over $\Spec\O_E$. It follows that
$N_0$ is the scheme theoretic closure of $N\otimes_{\O_E}E$ in $N$, i.e.\
$N_0=N^{\rm loc}$. Transferring this result to the standard model we
obtain the following result.

\begin{thm}\label{theorem5.7}
Assume that $(\vert\Gamma\vert, {\rm char}\, k)=1$ and that the conjecture
of De Concini-Procesi holds over $k$. Then the canonical flat model
$M^{\rm loc}$ can be described as the moduli scheme for the following
moduli functor on $({\rm Sch}/{\O}_E)$: The $S$-valued points are given by
${\O}_F\otimes_{{\cal O}_{F_0}} {\O}_S$-submodules $\F$ of
$\Lambda\otimes_{\O_{F_0}}{\O}_S$ which are locally direct summands as
${\O}_S$-modules, with $$
 {\rm{det}}(T\cdot{\rm Id}-\pi\ |\ {\cal F})\equiv
\prod_{\varphi}(T-\varphi(\pi))^{r_{\varphi}}
$$
 and
$$
\sum_{\sigma\in
\Gamma} F^t_h(\pi\ |\ \F)\cdot \wedge^t(\sigma Q_\phil (\pi\ |\ {\F}))=0,
$$
for all $\phil$ as before, and for $t+h=r-\sum_{\varphi,f(\varphi)\neq
0}r_{\varphi}+1$, $ t\geq 1, h\geq 0$.\endproof
\end{thm}

Let us fix integers $d,r$ and $e$ and let us state the following
conjecture.

\begin{conjecture}\label{conjecture5.8a}
Consider the closed subscheme $\overline N$ of ${\rm Mat}_{r\times r}$
over $k$ defined in (\ref{4.5}), $$ \overline N=\{ A\in {\rm Mat}_{r\times
r};\ A^e=0,\ {\rm det}(T\cdot I-A)\equiv T^r\}\ . $$ Then this scheme is
reduced.
\end{conjecture}

Note that when $r\leq e$, the first condition describing $\overline N$
follows from the second (Cayley-Hamilton) and in this case the statement
is a classical theorem of Kostant on the nilpotent variety. In a companion
paper to ours, J.~Weyman proves this conjecture when $e=2$, and for
arbitrary $e$ when char $k=0$. We note that the reduced subscheme of
$\overline N$ is simply the orbit closure corresponding to the partition
$(e^c, f)$ of $r$, where we have written as usual $r=c\cdot e+f$, $0\leq
f<e$. Now De Concini and Procesi [deC-P] have determined the length of the
intersection of any nilpotent orbit closure with the diagonal matrices (in
arbitrary characteristic). Therefore, if Conjecture \ref{conjecture5.8a}
holds we would obtain from their formula that
\begin{eq}\label{5.23}
{}\qquad{\rm dim}_k k[X_1,\ldots, X_r]/(e_1,\ldots, e_r, X_1^e,\ldots,
X_r^e)= {r!\over ((c+1)!)^f\cdot (c!)^{e-f}}\ .
\end{eq}
Here $e_1,\ldots, e_r$ are the elementary symmetric functions in the
indeterminates $X_1,\ldots,X_r$.

\begin{cor}\label{cor5.9b}
Let ${\bf r}= {\bf r}_{\rm min}(r,e)$. The corresponding standard model
for $GL_d$ is flat over $\Spec\O_E$ if Conjecture \ref{conjecture5.8a}
holds true.
\end{cor}

\begin{Proof}
If the conjecture holds true, we conclude from Theorem \ref{phismooth}
that $M\otimes_{\O_E}k$ is reduced. On the other hand, by Proposition
\ref{prop3.2} the generic and the special fiber of $M$ are irreducible of
the same dimension. The flatness of $M$ follows as in the proof of
Proposition \ref{normal} from EGA IV.3.4.6.1.
\endproof
\end{Proof}

\begin{cor}\label{cor5.8}
Assume that $e=2$ and order the embeddings so that $r_1\geq r_2$. Then the
canonical flat model $M^{\rm loc}$ represents the following moduli problem
on $(Sch/\O_E)$: The $S$-valued points are given by
$\O_F\otimes_{F_0}\O_S$-submodules ${\cal F}$ of
$\Lambda\otimes_{\O_{F_0}}\O_S$ which are locally direct summands as
$\O_S$-modules with $$ {\rm det}(T\cdot {\rm Id}-\pi\vert{\cal F})\equiv
(T-\varphi_1(\pi))^{r_1} (T-\varphi_2(\pi))^{r_2} $$ and $$
\wedge^{r_2+1}(\pi - \varphi_1(\pi)\cdot {\rm Id}\vert {\cal F})=0\ \
\hbox{if}\ r_1>r_2 \ .$$ Here in the last case we used $\varphi_1$ to
identify $F$ with $E$.
\end{cor}

\begin{Proof}
In this case $E$ is a Galois extension, namely $E=F_0$ if $r_1=r_2$ and
$E=F$ via $\varphi_1$ if $r_1>r_2$. Let us discuss the case when
$r_1>r_2$. In this case $\Gamma$ is trivial, and the second identity above
is the one in Theorem \ref{theorem5.7} corresponding to $f$ with $f(1)=1$,
$f(2)=0$ and to $t=r_2+1$, $h=0$. It follows that the conditions above
define a closed subscheme $M'$ of $M$ and it is easy to see that the
generic fibers coincide. The special fiber of $M'$ is defined by the
conditions
\begin{eq}\label{5.23}
{\rm det}(T\cdot {\rm Id}-\pi\vert {\cal F})\equiv T^r\ \ ,\ \
\wedge^{r_2+1}(\pi\vert {\cal F})=0\ \ .
\end{eq}
Consider the closed subscheme $N'$ of $N$ defined by the condition
$\wedge^{r_2+1}A=0$. Then the special fiber of $N'$ is given as
\begin{eq}\label{5.24}
\overline N'=\{A\in {\rm Mat}_{r\times r};\ A^2=0, \wedge^{r_2+1}A=0, {\rm
det}(T\cdot I-A)\equiv T^r\} .
\end{eq}
We then obtain a diagram of morphisms
\begin{eq}\label{5.25}
M'\buildrel\pi'\over\longleftarrow \widetilde
M'\buildrel\phi'\over\longrightarrow N'
\end{eq}
in which $\pi'$ is a $GL_r$-torsor and $\phi'$ is smooth of relative
dimension $rd$. According to Strickland ([St]), the special fiber of $N'$
is irreducible and reduced and in fact equal to the reduced closure of the
nilpotent orbit corresponding to the partition ${\bf s}=(2^{r_2},
1^{r_1-r_2})$, hence of dimension $r^2-(r_1^2+r_2^2)$. It follows that the
special fiber of $M'$ is irreducible and reduced of dimension
$dr-(r_1^2+r_2^2)={\rm dim}\ M'\otimes_{\O_E}E$. Hence $M'$ is $\O_E$-flat
and therefore coincides with $M^{\rm loc}$.

\medskip\noindent
Now assume $r_1=r_2$. In this case we are asserting that the standard
model is flat over ${\cal O}_{F_0}$. This follows via Corollary
\ref{cor5.9b} from Weyman's result on Conjecture \ref{conjecture5.8a}
concerning $e=2$.
\endproof
\end{Proof}

\medskip\noindent
{\it Remark:} The second identity above first appeared in [P] with the
purpose of defining a flat local model for the unitary group corresponding
to a ramified quadratic extension of $\Q_p$. It was the starting point of
the present paper.

\begin{Remark}\label{remark5.9}
{\rm  Corollary \ref{cor5.8} illustrates the fact that the
identities in Theorem \ref{theorem5.7} are extremely redundant. It
is an open problem to find a shorter list of identities, with
coefficients in $\O_E$, in suitable tensor powers of $\pi\vert{\cal
F}$ which describe the canonical flat model. }
\end{Remark}

\section{Applications to the affine
Grassmannians} \label{section5A} \setcounter{equation}{0}

\medskip\noindent
In this section we spell out some consequences of the results of the
previous sections, as they pertain to the special fibers. For the first
application we take up again the notation of section \ref{grassmannian};
in particular, we denote for a dominant coweight ${\bf s}$ of $GL_d $, by
${\cal O}_{\bf s}$ the corresponding orbit of $\tilde{\cal G}$ on
$\widetilde{\rm Grass}_k$.

\begin{thm}\label{theorem5A1}
 The reduced closure $\overline{\cal O}_{\bf s}$
is normal with rational singularities. Its singularities  are in fact
smoothly equivalent to singularities occurring in nilpotent orbit closures
for a general linear group.
\end{thm}

\begin{Proof}
Let ${\bf s} =(s_1,\ldots, s_d)$. After translation by a scalar matrix we
may assume $s_d\geq 0$. Let $e$ be any integer $\geq s_1$, and put $r=\sum
s_i$. Then ${\bf s}\in{\cal S}^0(r,e,d)$. Hence ${\cal O}_{\bf s}$ may be
identified with the corresponding stratum $M_{\bf s}$ of the special fiber
$\overline M$ of any standard model $M(\Lambda, {\bf r})$, where
$[F:F_0]=e$ and $\Lambda ={\cal O}_F^d$ and $\sum_{j=1}^e r_j=r$ (all of
them have identical special fibers). And the closure $\overline{\cal
O}_{\bf s}$ of ${\cal O}_{\bf s}$ can be identified with the closure
$\overline M_{\bf s}$ of $M_{\bf s}$ in $\overline M$. From Theorem
\ref{phismooth} it follows that $\overline M_{\bf s}$ is smoothly
equivalent to the closure of the corresponding nilpotent orbit $N_{\bf s}$
of $\overline N$, which by Mehta - van der Kallen is normal with rational
singularities.
\end{Proof}\endproof

\begin{Remark}\label{remark5A2} {\rm
 Results of this type (normality of Schubert varieties) have
been shown in positive characteristic in the context of Kac-Moody algebras
by Mathieu [Mat]. However, it is not clear whether the Schubert varieties
he considers have the same scheme structure as the $\overline{\cal O}_{\bf
s}$ considered here (the corresponding statement is known in
characteristic zero, by the integrability result of Faltings, comp.\ [BL],
app.\ to section 7). It follows from the methods of G\"ortz [G] that, once
a statement of the kind of Theorem \ref{theorem5A1} is known, it follows
that {\it any} Schubert variety in {\it any} parahoric flag variety for
$GL_d(k((\Pi))\, )$ is normal (and much more, comp.\ [G1]). For a
completely different approach see Faltings's paper [F].}
\end{Remark}

In fact, we can be more precise than in Theorem \ref{theorem5A1}. Let
$(V_e, \Pi_e)$ be the standard vector space of dimension $e$ over $k$,
with the standard regular nilpotent endomorphism. For $d\geq 1$ let

\begin{eq}\label{5A1}(W, \Pi)= (V_e, \Pi_e)^d\ \
.\end{eq} Let $0\leq r\leq ed$ and consider the projective scheme
$X=X(r,e,d)$ over $k$ which represents the following functor on the
category of $k$-algebras. It associates to a $k$-algebra $R$ the set of
$R$-submodules,
\begin{eqnarray*}
\{ {\cal F}\subset W\otimes_kR;\ && {\cal F}\ \hbox{is locally on ${\rm
Spec}\ R$ a direct summand,}\\  &&{\cal F}\ \hbox{is $\Pi$-stable and}\\
 && {\rm det}(T-\Pi\vert{\cal F})\equiv T^r\}\ .
\end{eqnarray*}

By the end of section \ref{grassmannian}, $X$ is the special fiber of any
standard model $M(\Lambda, {\bf r})$ where $[F:F_0]=e$, where $\Lambda
={\cal O}_F^d$ and $r=\sum r_{\varphi}$. We obtain as a special fiber of
the diagram (\ref{4.4}) the diagram

\begin{eq}\label{5A2}
X\buildrel\pi\over\longleftarrow \tilde X\buildrel\phi\over
\longrightarrow \overline N\ \ ,
\end{eq}
in which $\pi$ is a $GL_r$-torsor and $\phi$ is smooth of relative
dimension $rd$. It is equivariant with respect to the action of the
product group $\overline{\cal G}\times GL_r$, where $\overline{\cal
G}={\cal R}_{k[\Pi]/(\Pi^e)/k}(GL_d)$, cf.\ (\ref{3.8}). We therefore
obtain an extremely close relationship between affine Grassmannians and
nilpotent varieties. Indeed, $X$ is the Schubert variety in the affine
Grassmannian for $GL_d$ corresponding to the coweight $(e^c,
f,0,\ldots,0)$ (where we have written as usual $r=c\cdot e+f$, $0\leq
f<e$), and the image of $\tilde X$ is the closure of the nilpotent orbit
corresponding to the partition $(e^c,f)$ of $r$.

\begin{Remark}\label{Remark5A3}{\rm
For this result the fact that the nilpotent endomorphism $\Pi$ in
(\ref{5A1}) is ``homogeneous'' is essential. Indeed, let $r=e\geq 2$ and
change momentarily notations to consider the following inhomogeneous
example. Let

\begin{eq}\label{5A3}
(W, \Pi)= k[\Pi]/(\Pi^{e})\oplus k[\Pi]/(\Pi)\ \ ,
\end{eq}
and define $X$ as above. Then it is easy to see that any ${\cal F}\in X$
satisfies $${\rm span}\{ \Pi v_1,\ldots, \Pi^{e-1}v_1\}
\mathop{\subset}\limits_{\neq} {\cal F}\mathop{\subset}\limits_{\neq} W\ \
,$$ where $v_1$ resp.\ $v_2$ denotes the generator as $k[\Pi]$-module of
the first resp.\ second summand of (\ref{5A3}). Hence $X\simeq {\bf P}^1$.
Let $o\in X$ be the special point corresponding to $${\cal F}_o={\rm
span}\{ v_2, \Pi v_1,\ldots, \Pi^{e-1}v_1\}\ \ .$$ For ${\cal F}_o$, the
Jordan type of $\Pi\vert{\cal F}_o$ is $(e-1, 1)$ whereas for ${\cal
F}\neq {\cal F}_0$, the Jordan type of $\Pi\vert {\cal F}$ is $(e)$. It
follows that the fiber of $\phi$ through a point of $\pi^{-1}({\cal F}_o)$
has dimension equal to ${\rm dim}\ GL_r - {\rm dim}\ N_{(e-1,1)}=e+2$. On
the other hand, the fiber of $\phi$ through a point of
$\pi^{-1}(X\setminus \{ o\})$ has dimension equal to $({\rm dim}\ GL_r+1)
-{\rm dim}\ N_{(e)} = e+1$. Hence $\phi$ is not smooth in this case. }
\end{Remark}

\begin{Remark}\label{Remark5A3a}{\rm
Let us fix ${\bf r}$ with $r=\sum r_{\varphi}$, and let us consider a
standard model $M(\Lambda, {\bf r})$ with special fiber $X=X(r,e,d)$. Let
us order ${\bf r}=(r_1,\ldots, r_e)$. After extension of scalars from $k$
to $k'$ we obtain as special fiber of the diagram (\ref{5.8c}) the
following diagram with cartesian squares,
\begin{eq}\label{5A3a}
\matrix{ \overline{\cal M} & \longleftarrow & \overline{\tilde{\cal M}} &
\longrightarrow & \overline{\cal N} \cr \big\downarrow && \big\downarrow
&& \big\downarrow \cr X\otimes_kk' & \longleftarrow & \tilde X\otimes_kk'
& \longleftarrow & \overline N\otimes_kk' & . \cr}
\end{eq}
Here $\overline{\cal N}$ is the Springer resolution of the nilpotent orbit
closure corresponding to ${\bf t}={\bf r}^{\vee}$, comp.\ beginning of
section \ref{flat}. On the other hand, as Ng\^o pointed out to us, the
variety $\overline{\cal M}$ is an object which is well-known in the theory
of the affine Grassmannians, comp.\ [N-P]. Namely, let us introduce the
$e$ minuscule coweights $\mu_i=(1^{r_i}, 0^{d-r_i})$ of $GL_d$.
Corresponding to $\mu_i$ we have the Schubert variety $\overline{\cal
O}_{\mu_i}$ in the affine Grassmannian for $GL_d$ {\it over} $k'$. Then we
may identify the variety $\overline{\cal M}$ with the {\it convolution} in
the sense of Lusztig, Ginzburg, Mirkovic and Vilonen
\begin{eq}\label{5A3b}
\overline{\cal O}_{(\mu_1,\ldots, \mu_e)}:= \overline{\cal
O}_{\mu_1}\tilde\times\ldots \tilde\times \overline{\cal O}_{\mu_e}
\end{eq}
and the morphism from $\overline{\cal M}$ to $X\otimes_kk'$ factors
through a proper surjective morphism which may be identified with the
natural morphism ([N-P], \S 9)
\begin{eq}\label{5A3c}
m_{(\mu_1,\ldots, \mu_e)}:\overline{\cal O}_{(\mu_1,\ldots, \mu_e)}
\longrightarrow \overline{\cal O}_{\mu_1+\ldots +\mu_e}\ \ .
\end{eq}
Note that $\overline{\cal O}_{\mu_1+\ldots +\mu_e}$ is just the Schubert
variety corresponding to the coweight ${\bf t}={\bf r}^{\vee}$.
 }
\end{Remark}

We note the following consequence of (\ref{5A2}).

\begin{prop}\label{prop5A3}
If $r\leq e$, the scheme $X(r,e,d)$ is reduced and locally a complete
intersection. If $r>e$ and Conjecture \ref{conjecture5.8a} is true, then
$X$ is still reduced.
\end{prop}

\begin{Proof}
We argue with the special fiber $\overline M$ of a standard model as
described before. If $r\leq e$, then we may take ${\bf r}$ such that
$r_j\leq 1$, $\forall j$. The result then follows from Corollary
\ref{cor4.6} If $r>e$, then $X$ is smoothly equivalent to (an open
subscheme of) $\overline N$, and the assertion follows from a positive
answer to Conjecture \ref{conjecture5.8a}.
\end{Proof}
\endproof

\medskip
It is well-known that the Grassmannian over $k$ associated to $GL_d$ is
not reduced (this happens already for $d=1$). Recall ([BL]) that
\begin{eq}\label{5A4}
\widetilde{\rm Grass}_k= GL_d(k((\Pi))\, )/GL_d(k[[\Pi]])\ \ ,\end{eq}

where $GL_d(k((\Pi))\, )$ resp.\ $GL_d(k[[\Pi]])$ is the ind-group scheme
resp.\ group scheme which to a $k$-algebra $R$ associates $GL_d(R((\Pi))\,
)$ resp.\ $GL_d(R[[\Pi]])$, and where the quotient is taken in the
category of $k$-spaces and turns out to be an ind-scheme, [BL], 2.2. On
the other hand, the analogous quotient for $SL_d$ instead of $GL_d$ is an
ind-scheme which is reduced and even integral ([BL], 6.4.),
\begin{eq}\label{5A5}
\widetilde{\rm Grass}_k^{(0)}=SL_d(k((\Pi))\, )/SL_d(k[[\Pi]])\ \
.\end{eq} (At this point the blanket assumption in loc.\ cit.\ that char
$k=0$ is not used.) This means that $\widetilde{\rm Grass}_k^{(0)}$ can be
obtained as an increasing union of integral $k$-schemes. However, the
rather indirect proof of this fact in loc.cit.\ does not give an explicit
presentation of $\widetilde{\rm Grass}_k^{(0)}$ as such an increasing
union. Based on Proposition \ref{prop5A3} we are able to give such a
presentation, provided Conjecture \ref{conjecture5.8a} holds true. For
this recall ([BL], 2.3) that $\widetilde{\rm Grass}_k^{(0)}$ represents
the functor on $k$-algebras which to a $k$-algebra $R$ associates the set
of special lattices in $R((\Pi))^d$. (A {\it lattice} is a
$R[[\Pi]]$-submodule $W$ of $R((\Pi))^d$ with $\Pi^fR[[\Pi]]^d\subset
W\subset \Pi^{-f}R[[\Pi]]^d$ for some $f$ and such that the $R$-module
$\Pi^{-f}R[[\Pi]]^d /W$ is projective. A lattice $W$ is {\it special} if
the lattice $\wedge^dW$ in $\wedge^dR((\Pi))^d=R((\Pi))$ is trivial, i.e.\
equal to $R[[\Pi]]$). Let $\widetilde{\rm Grass}_k[0]$ be the space of
lattices of total degree 0 (i.e.\ the rank of $\Pi^{-f}R[[\Pi]]^d/W$ is
equal to $fd$. Then $\widetilde{\rm Grass}_k[0]$ is a connected component
of $\widetilde{\rm Grass}_k$ and
\begin{eq}\label{5A6}
\widetilde{\rm Grass}_k^{(0)} =(\widetilde{\rm Grass}_k[0])_{\rm red}\ \ ,
\end{eq}
 (cf.\
[BL], 2.2.\ and 6.4.).

Let $\tilde X_f$ be the subscheme of $\widetilde{\rm Grass}_k[0]$ which
parametrizes the lattices $W$ with $\Pi^fR[[\Pi]]^d\subset W\subset
\Pi^{-f}R[[\Pi]]^d$ and let $X_f$ be the closed subscheme of $W$ in
$\tilde X_f$ such that ${\rm det}(T-\Pi\vert (W/\Pi^fR[[\Pi]]^d))\equiv
T^{fd}$. We have an exact sequence $$0\to\Pi^fR[[\Pi]]^d /
\Pi^{f+1}R[[\Pi]]^d\to W/\Pi^{f+1}R[[\Pi]]^d\to W/\Pi^fR[[\Pi]]^d\to 0\
,$$ which implies that $${\rm det}(T-\Pi\vert (W/\Pi^{f+1} R[[\Pi]]^d))
={\rm det} (T-\Pi\vert (W/\Pi^fR[[\Pi]]^d))\cdot T^d\ \ .$$ We therefore
obtain a chain of closed embeddings of $k$-schemes
\begin{eq}\label{5A7}
X_0\subset X_1\subset\ldots
\end{eq}

\begin{prop}\label{prop5A4}
 For $d\leq 2$, the chain (\ref{5A7}) presents $\widetilde{\rm
Grass}_k^{(0)}$ as an increasing union of integral $k$-schemes. The same
is true for arbitrary $d$, if Conjecture \ref{conjecture5.8a} holds true.
In particular, this holds (by the theorem of Weyman) if char $k=0$.
\end{prop}

\begin{Proof}
We note that $X_f=X(fd, 2f, d)$, hence is integral if $d\leq 2$ and for
arbitrary $d$, if Conjecture \ref{conjecture5.8a} holds true. Hence
$X_f=(\tilde X_f)_{\rm red}$. The claim follows from $\widetilde{\rm
Grass}_k[0] =\lim\limits_{\longrightarrow}\tilde X_f$ and $$\widetilde{\rm
Grass}_k^{(0)}=\lim\limits_{\longrightarrow} (\tilde X_f)_{\rm red}
=\lim\limits_{\longrightarrow} X_f\ \ .\eqno\endproof$$
\end{Proof}

\section{The complex of nearby cycles} \label{section6}
\setcounter{equation}{0}

In this section we suppose that the residue field $k$ of $\O_E$ is finite
and shall aim for an expression for the complex of nearby cycles of a
standard model for $GL_d$ corresponding to ${\bf r}$, by exploiting in
more depth the resolution of singularities given by the scheme $\N$ in
(\ref{5.3}).

Recall that $K$ denotes the Galois hull of $F$ in $F_0^{\rm sep}$ and that
$k'$ is the residue field of $\O_K$. Also, $\Gamma$ denotes the Galois
group of $K$ over $E$. We fix a prime number $\ell$ which is invertible in
$\O_E$ and denote by $R\psi =R\psi\overline{\Q}_{\ell}$ the complex of
nearby cycles of the ${\O_E}$-scheme $M$. Since $M^{\rm loc}$ is the
scheme-theoretic closure of the generic fiber in $M$, this complex has
support in $\overline M^{\rm loc}:=M^{\rm loc}\otimes_{\O_E}k$. The
complex is equipped with an action of ${\rm Gal}( F_0^{\rm sep}/E)$ which
lifts the action on $\overline M^{\rm loc}$. We also fix a square root of
the cardinality $\vert k\vert$ in $\overline{\Q}_{\ell}$.

\begin{thm}\label{theorem6.1}
There is an isomorphism between perverse sheaves pure of weight zero,
$$
R\psi [{\rm dim}\ \overline M^{\rm loc}] \hbox{$\left( {1\over 2} {\rm
dim}\ \overline M^{\rm loc}\right)$} =\bigoplus_{{\bf s}\leq {\bf t}}{\cal
M} _{\bf s}\otimes IC_{M_{\bf s}}\ ,
$$
 where $IC_{M_{\bf s}}$ denotes the
intermediate extension of the constant sheaf\break
 $\overline{\Q}_{\ell} [{\rm dim}\ M_{\bf s}] \left( {1\over 2} {\rm
dim}\ M_{\bf s}\right)$, equipped with the action of ${\rm Gal}( F_0^{\rm
sep} /E)$ which factors through ${\rm Gal}(\overline k/k)$, and ${\cal
M}_{\bf s}$ is a $\overline{\Q}_{\ell}$-vector space equipped with an
action of ${\rm Gal}( F_0^{\rm sep}/E)$.  The action of ${\rm Gal}(
F_0^{\rm sep}/E)$ on ${\cal M}_{\bf s}$ factors through $\Gamma$, and the
degree of this representation is given as a Kostka number, $$ m_{\bf
s}={\rm dim}\ {\cal M}_{\bf s}= K_{{\bf s}^{\vee}, {\bf r}} $$ (cf.\ [Mc],
p.\ 115).
\end{thm}

By making use of naturality properties of the complex of nearby cycles,
the diagram in Theorem \ref{mcan} (iii) with smooth horizontal morphisms
reduces us to proving the corresponding statements for $N^{\rm loc}$
instead of $M^{\rm loc}$ (or $M$). To be more precise, let us change
notations and denote now by $R\psi$ the complex of nearby cycles of the
$\O_E$-scheme $N^{\rm loc}$. We wish to prove the formula

\begin{eq}\label{6.1}
R\psi [{\rm dim}\ \overline N^{\rm loc}] \hbox{$\left( {1\over 2} {\rm
dim}\overline N^{\rm loc}\right) $}=\bigoplus_{{\bf s}\leq {\bf t}} {\cal
M} _{\bf s}\otimes IC_{N_{\bf s}}\ \ ,
\end{eq}
with ${\cal M}_{\bf s}$ as above.

The left hand side is a perverse sheaf of weight zero on $\overline N^{\rm
loc}$ ([I] Thm.\ 4.2 and Cor.\ 4.5), which is $GL_r$-equivariant. By
[BBD], 5.3.8 its inverse image $R\overline\psi$ on $\overline N^{\rm
loc}\otimes_k\overline k$ is semisimple and its simple constituents are
all of the form $IC_{N_{\bf s}\otimes_k\overline k}$, for some ${\bf
s}\leq {\bf t}$, since $N_{\bf s}\otimes_k\overline k$ admits no
non-trivial $GL_r$-equivariant irreducible local system. Therefore, the
isotypical decomposition of $R\overline\psi$ has the form

\begin{eq}\label{6.1a}
R\overline\psi=\bigoplus_{{\bf s}\leq{\bf t}}K_{\bf s}\ \ ,
\end{eq}

\noindent where $K_{\bf s}$ is a multiple of $IC_{N_{\bf s}\otimes_k\overline k}$.
Obviously $K_{\bf s}$ is of the form

\begin{eq}\label{6.1b}
K_{\bf s}={\cal M}_{\bf s}\otimes IC_{N_{\bf s}\otimes_k\overline k}\ \ ,
\end{eq}

\noindent where ${\cal M}_{\bf s}$ is a $\overline{\Q}_{\ell}$-vector space with an
action of the inertia group $I\subset {\rm Gal}( F_0^{\rm sep}/E)$. The
fact that both sides of (\ref{6.1b}) come by extension of scalars from
$k$, implies now that ${\rm Gal}( F_0^{\rm sep}/E)$ acts on ${\cal M}_{\bf
s}$ and that we have a decomposition of the form (\ref{6.1}).

It remains to show that the restriction of ${\cal M}_{\bf s}$ to ${\rm
Gal}( F_0^{\rm sep}/K)$ is trivial and to determine its degree.  By
Deligne [D], Prop.3.7., we have

\begin{eq}\label{6.2}
R\psi \otimes_kk'=R\psi'\ \
\end{eq}

\noindent with $R\psi'$ the complex of nearby cycles of the $\O_K$-scheme
$N^{\rm loc}\otimes_{\O_E}\O_K$. Now $N^{\rm loc}\otimes_{\O_E}K$ has the
smooth model $\N$ which maps via $\mu$ to $N^{\rm loc}\otimes_{\O_E}\O_K$.
Since the complex of nearby cycles of a smooth scheme is the constant
sheaf placed in degree 0, the functoriality with respect to push-forward
under a proper morphism gives a natural identification of complexes on
$\overline N^{\rm loc}\otimes_kk'$,

\begin{eq}\label{6.3}
R\psi'= R\overline{\mu}_*\overline{\Q}_{\ell}\ \ .
\end{eq}

In particular, the action of ${\rm Gal}(F_0^{\rm sep}/K)$ on $R\psi'$ is
through ${\rm Gal}(\overline k/k')$. Now the morphism
$\overline\mu:\overline{\cal N}\otimes_kk'\to \overline N^{\rm
loc}\otimes_kk'$ comes by base change from the moment map
$\overline{\overline\mu}$ for the variety $\overline{\hbox{\boldmath $\cal
F$}}$ of partial flags over $k$, cf.\ (\ref{5.2}). This map is semi-small
with all strata $N_{\bf s}$ of $\overline N^{\rm loc}$ relevant, comp.\
[BM]. Base changing the decomposition (\ref{6.1}) from $k$ to $k'$ and
identifying the left hand side with $R
{\overline\mu}_*\overline{\Q}_{\ell}[{\rm dim}\overline N^{\rm
loc}]({1\over 2}{\rm dim}\ \overline N^{\rm loc})$, we see that

\begin{eq}\label{6.3a}
{\cal M}_{\bf s}=R^{2d_{\bf s}}\overline\mu_*\overline{\Q}_{\ell}(d_{\bf
s})_{\vert N_{\bf s}\otimes_k\overline k}\ \ ,
\end{eq}

\noindent as representations of ${\rm Gal}(F_0^{\rm sep}/K)$, i.e.\ of
${\rm Gal}(\overline k/k')$. Here $d_{\bf s}={1\over 2}{\rm codim}\ N_{\bf
s}$ is the relative dimension of $\overline\mu$ over $N_{\bf s}$. But by
Spaltenstein [Sp] all irreducible components of
$\overline{\overline\mu}^{-1}(N_{\bf s}\otimes_k\overline k)$ are defined
over $k$, hence ${\rm Gal}(\overline k/k')$ acts trivially on the right
hand side of (\ref{6.3a}). This shows that ${\rm Gal}(F_0^{\rm sep}/K)$
acts trivially on ${\cal M}_{\bf s}$. (Spaltenstein works over an
algebraically closed field, but his results are valid over an arbitrary
field, comp.\ [HS], \S 2.) A different argument, pointed out by G.Laumon,
is to transpose the result of Braverman and Gaitsgory [BG] from ${\bf C}$
to a finite field and appeal to [HS], Corollary 2.3.

For the degree $m_{\bf s}$ of ${\cal M}_{\bf s}$ Borho and MacPherson give
the formula [BM], 3.5.,

\begin{eq}\label{6.5}
m_{\bf s}={\rm dim}\ {\rm Hom}_{S_r}(\chi^{\bf s}, {\rm Ind}^{S_r}_{S_{\bf
r}}({\rm sgn}))\ .
\end{eq}

\noindent Here $S_{\bf r}=S_{r_1}\times \ldots \times S_{r_e}$ is a
subgroup of the symmetric group $S_r$ and ${\rm sgn}$ is the sign
character on all factors. Furthermore, $\chi^{\bf s}$ denotes the unique
irreducible representation of $S_r$ which occurs both in ${\rm
Ind}^{S_r}_{S_{\bf s}}(1)$ and in ${\rm Ind}^{S_r}_{S_{\bf s}}({\rm
sgn})$, cf.\ [Mc], p.\ 115. By loc.cit.\ (\ref{6.5}) can be identified
with the Kostka number occurring in Theorem \ref{theorem6.1}. Note that
this agrees with the formula of Braverman and Gaitsgory [BG], Cor.\ 1.5.,

\begin{eq}\label{6.6}
m_{\bf s}={\rm dim}\ V({\bf s}^{\vee})_{\bf r}\ \ ,
\end{eq}

\noindent (where, however, in their formula $V({\bf P})_{\bf d}$ should be
replaced by $V({\bf P}^{\vee})_{\bf d}$).

Here $V({\bf s}^{\vee})$ denotes the rational representation of $GL_e$,
$$V({\bf s}^{\vee})= {\rm Hom}_{S_r} (\chi^{{\bf s}^{\vee}}, ({\rm
nat}_{GL_e})^{\otimes r})$$ where ${\rm nat}$ is the natural
representation of $GL_e$ and $V({\bf s}^{\vee})_{\bf r}$ the weight space
corresponding to ${\bf r}$. By [Mc], p.\ 163 the character of $V({\bf
s}^{\vee})$ is given by the Schur function $s_{{\bf s}^{\vee}}$ and hence
(\ref{6.6}) is indeed equal to $K_{{\bf s}^{\vee}, {\bf r}}$ by [Mc], p.\
101. (Both sources [BM] and [BG] work over ${\bf C}$, but can be
transposed to the present context.)\endproof

\begin{Remark}\label{remark6.1a}{\rm
We have used here Theorem \ref{phismooth} in order to deduce Theorem
\ref{theorem6.1} from the smoothness of ${\cal N}$ and the fact that the
Springer resolution of the nilpotent orbit closure corresponding to ${\bf
t}={\bf r}^{\vee}$ is semi-small. Of course, Theorem \ref{phismooth}
implies also that the scheme ${\cal M}$ is smooth and that the special
fiber of ${\cal M}$ is a semi-small resolution of $\overline M_{\bf t}$.
On the other hand, this last fact has a direct proof, cf.\ [N-P], Lemma
9.3. In fact, the Remark \ref{Remark5A3a} establishes an equivalence
between the semi-smallness properties.
 }
\end{Remark}

\begin{Remark}\label{ngoremark}
{\rm (B. C. Ng\^o) Denote by $R\psi'$ the complex of nearby cycles
of the $\O_K$-scheme $M^{\rm loc}\otimes_{\O_E}\O_K$. Let
us fix an $\O_F$-basis of the module $\Lambda$.
Recall that then by (\ref{5A3b}), the special fiber $\overline{\cal M}$
can be identified with the convolution
$$
\overline{\cal
O}_{\mu_1}\tilde\times\ldots \tilde\times \overline{\cal O}_{\mu_e}
$$
and the morphism $\overline{\cal M}\to \overline{M}^{\rm loc}\times_k k'$
with the natural morphism
$$
m_{(\mu_1,\ldots, \mu_e)}:\overline{\cal
O}_{\mu_1}\tilde\times\ldots \tilde\times \overline{\cal O}_{\mu_e}
\to \overline{\cal O}_{\mu_1+\cdots +\mu_e}
$$
of (\ref{5A3c}); here $\overline{\cal O}_{\mu}$
denotes the reduced closure of the orbit which corresponds to
the dominant coweight $\mu$ in the affine Grassmanian for
$GL_d$ over $k'$.
The same arguments as in the proof of Theorem
\ref{theorem6.1}, applied to the morphism $\pi: {\cal M}\to M^{\rm loc}\otimes_{\O_E}\O_K$
in place of ${\cal N}\to N^{\rm loc}\otimes_{\O_E}\O_K$, show now that
\begin{eq}\label{6.6a}
R\psi' [{\rm dim}\ \overline M^{\rm loc}] \hbox{$\left( {1\over 2} {\rm
dim}\ \overline M^{\rm loc}\right)$} = R\overline\pi_* \overline{\Q}_l
[{\rm dim}\ \overline {\cal M}] \hbox{$\left( {1\over 2} {\rm
dim}\ \overline {\cal M}\right)$}.
\end{eq}
By the above discussion, the right hand side of (\ref{6.6a}) is the
convolution
$
IC_{{\cal O}_{\mu_1}}*\cdots * IC_{{\cal O}_{\mu_e}}
$. Hence,
we obtain a relation between nearby cycles and convolution:
\begin{eq}\label{conv}
R\psi' [{\rm dim}\ \overline M^{\rm loc}] \hbox{$\left( {1\over 2} {\rm
dim}\ \overline M^{\rm loc}\right)$}=
IC_{{\cal O}_{\mu_1}}*\cdots * IC_{{\cal O}_{\mu_e}}.
\end{eq}

Recall that there is an equivalence of tensor categories
between the category of $\tilde {\cal G}_{k'}$-equivariant pure perverse
$\overline{\Q}_l$-sheaves of weight $0$ on the affine Grassmanian
$\widetilde {\rm Grass}_{k'}$ (with tensor structure
given by the convolution product)
and the category of finite dimensional
$\overline{\Q}_l$-representations of the Langlands dual group
$GL_d=\widehat{GL}_d$ (see [Gi], and especially [M-V] \S 7).
Under this equivalence the perverse sheaf $IC_{{\cal O}_\mu}$
corresponds to the representation $V(\mu)$ of
${GL}_d$ of highest weight
$\mu$, and
the convolution in (\ref{conv}) to
the tensor product
$$
V(\mu_1)\otimes \cdots \otimes V(\mu_e).
$$
This
decomposes
$$
V(\mu_1)\otimes \cdots \otimes V(\mu_e)\ =\bigoplus_{\lambda\leq \mu _1+\cdots+\mu_e}
{\cal M}_{\lambda}\otimes V(\lambda) \ ,
$$
where ${\cal M}_\lambda$ is a finite dimensional
$\overline\Q_l$-vector space.
Using again the above equivalence of categories and (\ref{conv}),
we obtain
\begin{eq} \label{ngo}
R\psi' [{\rm dim}\ \overline M^{\rm loc}] \hbox{$\left( {1\over 2} {\rm
dim}\ \overline M^{\rm loc}\right)$}=\bigoplus_{\lambda\leq \mu _1+\cdots+\mu_e}
{\cal M}_{\lambda}\otimes IC_{{\cal O}_\lambda} \ .
\end{eq}
The right hand side of (\ref{ngo}) corresponds to the expression
in Theorem \ref{theorem6.1}. Indeed, $\{\lambda \ |\ \lambda\leq \mu _1+\cdots+\mu_e\}$
corresponds to
$\{{\bf s}\ |\ {\bf s}\leq {\bf r}^\vee\}$, and we can see directly that
the Littlewood-Richardson number $\dim({\cal M}_\lambda)$ is equal to the Kostka
number $K_{{\bf s}^\vee, {\bf r}}$ of \ref{theorem6.1}.}
\end{Remark}

\begin{Remark}
{\rm As before, let us choose an ordering of the set of embeddings $\phi:
F\to F_0^{\rm sep}$. The Galois group $\Gamma={\rm Gal}(K/E)$ can be
identified with the subgroup of elements $\sigma$ of the symmetric group
$S_e$ which satisfy $r_{\sigma(i)}=r_{i}$. The group $\Gamma$ acts by
permutation of the factors on the tensor product $V(\mu_1)\otimes \cdots
\otimes V(\mu_e)$. Let us denote by $\rho$ the corresponding
representation $$ \rho\ :\ \Gamma\to {GL}(V(\mu_1)\otimes \cdots \otimes
V(\mu_e)). $$ Since the permutation action commutes with the action of
${GL}_d$ on the tensor product, the representation $\rho$ decomposes as $$
\rho=\bigoplus_{\lambda\leq \mu _1+\cdots+\mu_e} \rho_{\lambda}\otimes
{\rm id}_{V(\lambda)} $$ where $\rho_\lambda$ is a representation of
$\Gamma$ on the vector space ${\cal M}_\lambda$. We conjecture that the
representation of ${\rm Gal}(K/E)$ on ${\cal M}_{\bf s}$ (see Theorem
\ref{theorem6.1}) is isomorphic to  $\rho_\lambda$ (with $\lambda$
corresponding to $\bf s$).

After our choice of an $\O_F$-basis
of $\Lambda$, the generic fiber
of $M^{\rm loc}$ can be identified with
the product of Grassmanians
$$
M^{\rm loc}\otimes_{\O_K}K=
\prod_{i=1}^e{\rm Grass}_{r_i}(K^d).
$$
For $\sigma\in \Gamma$, permutation of the factors
gives an isomorphism of $K$-schemes
$$
\kappa(\sigma): M^{\rm loc}\otimes_{\O_E}K \to
M^{\rm loc}\otimes_{\O_E}K.
$$
The isomorphism $\kappa(\sigma)$ induces an automorphism
of the sheaf of vanishing cycles $R\psi'$.
Therefore, by (\ref{conv}) we obtain
a ``commutativity" isomorphism
$$
\kappa(\sigma): IC_{{\cal O}_{\mu_1}}*\cdots * IC_{{\cal O}_{\mu_e}}
\to IC_{{\cal O}_{\mu_1}}*\cdots * IC_{{\cal O}_{\mu_e}}=
IC_{{\cal O}_{\mu_{\sigma(1)}}}*\cdots * IC_{{\cal O}_{\mu_{\sigma(e)}}}.
$$
(notice here that, since $\sigma\in\Gamma$, $\mu_{\sigma(i)}=\mu_{i}$).
As was pointed out by Ng\^o, the conjecture follows, if $\kappa(\sigma)$
coincides with the isomorphism given using the permutation $\sigma$
and the ``commutativity constraint" (for the
tensor category of perverse sheaves of Remark \ref{ngoremark})
of [M-V].}
\end{Remark}

\begin{Remark}\label{remark6.2}
{\rm Let $x$ be a point of $M^{\rm loc}$ with values in a finite extension
${\bf F}_q$ of $k$. Let ${\rm Tr}^{ss}({\rm Fr}_q, R\psi_x^M)$ be the
semi-simple trace of the geometric Frobenius on the stalk at $x$ of the
complex of nearby cycles of $M^{\rm loc}$. It should be possible to
transfer the conjecture of Kottwitz [HN] to this case to obtain a
group-theoretic expression for this (the group $R_{F/F_0}(GL_d)$ relevant
for this conjecture here is not split so that the original formulation
does not apply directly). Note that

\begin{eq}\label{6.7}
{\rm Tr}^{ss}({\rm Fr}_q, R\psi_x^M)= {\rm Tr}^{ss}({\rm Fr}_q, R\psi_y^N)
\end{eq}
where $y=\phi(\tilde x)$ for an arbitrary point $\tilde x\in \widetilde
M^{\rm loc}({\bf F}_q)$ mapping to $x$, and where $R\psi^N$ is the complex
of nearby cycles of $N^{\rm loc}$. The expression (\ref{6.7}) only depends
on the stratum  $M_{\bf s}$ containing $x$ resp.\ $N_{\bf s}$ containing
$y$ and may therefore be denoted by ${\rm Tr}^{ss}(Fr_q, R\psi^N_{\bf
s})$.

We also consider the analogous semi-simple traces for points with values
in a finite extension ${\bf F}_q$ of $k'$,

\begin{eqnarray}\label{6.7a}
&&{\rm Tr}^{ss} (Fr_q, R\psi_x^{M\otimes K}),\ {\rm Tr}^{ss}(Fr_q,
R\psi_x^{N\otimes K})\ ,\\ &&{\rm Tr}^{ss}(Fr_q, R\psi_{\bf s}^{M\otimes
K}),\ {\rm Tr}^{ss}(Fr_q, R\psi_{\bf s}^{N\otimes K})\ .\nonumber
\end{eqnarray}
(these semi-simple traces may differ from the preceding ones when $K/E$ is
ramified).

>From Theorem \ref{theorem6.1} we obtain

\begin{eq}\label{6.8}
{\rm Tr}^{ss}({\rm Fr}_q;\ R\psi^{N\otimes K}_{\bf s}) = q^{{1\over 2}{\rm
dim}\ \overline N^{\rm loc}} \sum_{{\bf s}\leq {\bf s}'\leq {\bf t}}
m_{{\bf s}'}\cdot {\rm Tr}^{ss}({\rm Fr}_q, (IC_{N_{{\bf s}'}})_{\bf s})\
,
\end{eq}
provided ${\bf F}_q\supset k'$. By Lusztig [L], the entity ${\rm
Tr}^{ss}({\rm Fr}_q, (IC_{N_{{\bf s}'}})_{\bf s})$ can be expressed in
terms of Kostka polynomials $K_{{\bf s}', {\bf s}}(q^{-1})$, comp.\ also
[Mc], p.\ 245.

Consider the spectral sequence of nearby cycles,

\begin{eq}\label{6.9}
E_2^{pq}=H^p(M^{\rm loc}\otimes_k\overline k, R^q\psi^M)\Rightarrow
H^{p+q}(M^{\rm loc}\otimes_{\O_E}F_0^{\rm sep}, \overline{\Q}_{\ell})\ .
\end{eq}

\par\noindent
We obtain an identity of semi-simple traces (cf.\ [HN])

\begin{eq}\label{6.10}
{}\ \ {\rm Tr}^{ss}({\rm Fr}_q, H^*(M^{\rm loc}\otimes_{\O_E} F_0^{\rm
sep}, \overline{\Q}_{\ell}))= \sum_{x\in M^{\rm loc}({\bf F}_q)} {\rm
Tr}^{ss}({\rm Fr}_q, R\psi_x^M) .
\end{eq}

\par\noindent
Similarly we may consider the spectral sequence of nearby cycles for
$M^{\rm loc}\otimes_{{\cal O}_E}{\cal O}_K$. We then obtain an identity of
semi-simple traces (as ${\rm Gal}(F_0^{\rm sep}/K)$-modules), provided
that ${\bf F}_q\supset k'$,

\begin{eq}\label{6.10a}
{}\, \, \,\,\,\,\,\,\,\,{\rm Tr}^{ss}(Fr_q; H^*((M{\otimes}_{{\cal O}_E}
{\cal O}_K ){\otimes}_{{\cal O}_K} F_0^{\rm sep},
\overline{\Q}_{\ell})){=} {\mathop{\sum}\limits_{x\in M({\bf F}_q)}} {\rm
Tr}^{ss}({\rm Fr}_q, R\psi_x^{M\otimes K}).
\end{eq}

\par\noindent
However,

$$M^{\rm loc}\otimes_{\O_E}K= \prod_{\varphi}{\rm
Grass}_{r_{\varphi}}(V_{\varphi})\ \ ,$$ cf.\ (\ref{2.7}), and hence the
left hand side of (\ref{6.10a}) is known. Taking into account (\ref{6.7})
we obtain therefore from (\ref{6.10}) a combinatorial identity involving
Kostka polynomials, Kostka numbers etc. It might be interesting to
identify this combinatorial identity explicitly.
 }
\end{Remark}

\section{Local models of $EL$-type}\label{section7}
\setcounter{equation}{0}

In this section we use the following notation (following [RZ]).
\smallskip

$F_0$ a complete discretely valued field with perfect residue field

$F$ a finite separable field extension of $F_0$,

$B$ a simple algebra with center $F$,

$V$ a finite $B$-module,

$G=GL_B(V)$, as algebraic group over $F_0$,

$\mu: {\bf G}_{mK}\to G_K$ a one parameter subgroup, defined over some
sufficiently big extension $K$ contained in a fixed separable closure
$F_0^{\rm sep}$ of $F_0$, given up to conjugation. We assume that the
eigenspace decomposition of $V\otimes_{{\bf Q}_p}K$ is given by

\begin{eq}\label{7.1}
V\otimes_{F_0}K= V_0\oplus V_1\ \ .
\end{eq}

$E$ the field of definition of the conjugacy class of $\mu$.

${\cal L}$ a periodic $O_B$-lattice chain in $V$.
\smallskip

\par\noindent In the sequel we denote by $\O_{F_0}, \O_F, \O_B, \O_E$ the respective
rings of integers.

\medskip\noindent To these data we associate the following functor on
$(Sch/\O_E)$: The $S$-valued points are given by
\begin{enumerate}
\item[i)] a functor $\Lambda\mapsto t_{\Lambda}$ to the category of
$\O_B\otimes_{\O_{F_0}}\O_S$-modules on $S$
\item[ii)] a morphism of functors
$\varphi_{\Lambda}:\Lambda\otimes_{\O_{F_0}}\O_S\to t_{\Lambda}$.
\end{enumerate}
\par\noindent
The requirements on these data are:

\medskip\noindent a) $t_{\Lambda}$ is locally on $S$ a free $\O_S$-module
of finite rank, and we have the following identity of polynomial functions
on $\O_B$:

\begin{eq}\label{7.2}
{\rm det}_{\O_S}(a\vert t_{\Lambda})={\rm det}_K(a\vert V_1)\ \ .
\end{eq}
b) $\varphi_{\Lambda}$ is surjective, for all $\Lambda\in {\cal L}$.

\medskip\noindent
We remark that the standard models considered in sections \ref{GLd} --
\ref{section6} are a special case. Indeed, let $B=F$ and let ${\cal L}$
consist of the $F^{\times}$-multiples of a fixed $\O_F$-lattice
$\Lambda_0$ in the $d$-dimensional $F$-vector space $V$ and assume that
under the decomposition (\ref{2.5}) of $F\otimes_{F_0}F_0^{\rm sep}$ we
have

\begin{eq}\label{7.3}
V\otimes_{F_0}F_0^{\rm sep}{=}\bigoplus_{\varphi}V_{\varphi}  ,\
V_0\otimes_KF_0^{\rm sep}{=}\bigoplus_{\varphi} V_{0,\varphi} ,\
V_1\otimes_KF_0^{\rm sep}{=}\bigoplus_{\varphi}V_{1,\varphi}  ,
\end{eq}
with ${\rm dim}_{F_0^{\rm sep}}V_{0,\varphi}=r_{\varphi}$. If
$(t_{\Lambda}, \varphi_{\Lambda})_{\Lambda}$ is an $S$-valued point of the
moduli problem above, then ${\cal F}={\rm Ker}\ \varphi_{\Lambda_0}$ is an
$S$-valued point of the standard model for $GL_d$ corresponding to ${\bf
r}=(r_{\varphi})$ and this establishes an isomorphism of moduli problems.
Furthermore, $E$ has indeed the description given in section \ref{GLd} in
terms of ${\bf r}$.

In general, the above moduli problem is representable by a projective
scheme over $\Spec\O_E$ which is a closed subscheme of a form over
$\O_E$ of a product of Grassmannians. Let us make this more precise.

Let us first consider the case where $B=F$ is a totally ramified extension
of degree $e$ of $F_0$. Let $v_1,\ldots, v_d$ be a basis of $V$. For
$i=0,\ldots, d-1$ consider the $\O_F$-lattice $\Lambda_i$ in $V$ spanned
by

\begin{eq}\label{7.4}
\Lambda_i={\rm span}\ \{ \pi^{-1}v_1,\ldots, \pi^{-1} v_i, v_{i+1},\ldots,
v_d\}\ \ .
\end{eq}
Here $\pi$ denotes a uniformizer in $F$. This yields a complete periodic
lattice chain

\begin{eq}\label{7.5}
\ldots\to \Lambda_0\to \Lambda_1\to\ldots\to \Lambda_{d-1}\to
\pi^{-1}\Lambda_0\to\ldots
\end{eq}

Choose $I=\{ i_0<i_1<\ldots < i_{\ell}\}\subset \{ 0,\ldots, d-1\}$. Then
the periodic lattice ${\cal L}$ is isomorphic to the subchain of
(\ref{7.5}) where only $\Lambda_i$ with $i\in I$ are kept, for suitable
$I$. The points of the above functor with values in a $\O_E$-scheme $S$
can now be interpreted as the isomorphism classes of commutative diagrams

$$
\matrix{ \Lambda_{i_0, S} & \longrightarrow & \Lambda_{i_1,S} &
\longrightarrow & \ldots & \longrightarrow & \Lambda_{i_{\ell}, S} &
\longrightarrow & \pi^{-1}\Lambda_{i_0, S} \cr \big\uparrow &&\big\uparrow
&&&&\big\uparrow&&\big\uparrow \cr {\cal F}_0 & \longrightarrow & {\cal
F}_1 & \longrightarrow & \ldots & \longrightarrow & {\cal F}_{\ell}
&\longrightarrow & \pi^{-1}{\cal F}_0 & , \cr}
$$
where $\Lambda_{i_j,S}$
is $\Lambda_{i_j}\otimes_{\O_{F_0}}\O_S$ and where the ${\cal F}_j$ are
$\O_F\otimes_{\O_{F_0}}\O_S$-submodules which locally on $S$ are direct
summands of $\Lambda_{i_j,S}$ as $\O_S$-modules and where we have an
identity of polynomial functions on $\O_F$,

\begin{eq}\label{7.6}
{\rm det}(a\vert {\cal F}_j)={\rm det}_K(a\vert V_0)\ \ ,\ \
j=0,\ldots,\ell\ \ .
\end{eq}
Suppose that $K$ is a sufficiently big Galois extension so that

\begin{eq}\label{7.7}
F\otimes_{F_0}K= \bigoplus_{\varphi: F\to K}K\ \ .
\end{eq}
We obtain corresponding decompositions

\begin{eq}\label{7.8}
V\otimes_{F_0}K=\bigoplus_{\varphi}V_{\varphi}\ ,\ V_0\otimes_{F_0}K=
\bigoplus_{\varphi}V_{0,\varphi}\ .
\end{eq}
Then ${\rm dim}_KV_{\varphi}=d$ for all $\varphi$. Let $r_{\varphi}={\rm
dim}_KV_{0,\varphi}$. Then the determinant condition (\ref{7.6}) can be
rewritten as

\begin{eq}\label{7.9}
{\rm det}(a\vert{\cal F}_j)= \prod_{\varphi}\varphi(a)^{r_{\varphi}}\ \ .
\end{eq}
In other words, each ${\cal F}_i$ is an $S$-valued point of the standard
model for $GL_d$ corresponding to ${\bf r}=(r_{\varphi})$ (and
$\Lambda_i$). Our functor is representable by a closed subscheme of the
product over $\O_E$ of these standard  models, one for each
$i=0,\ldots,\ell$.

Now let us consider the general case. Let $\Fsmile_0$ be the completion of
the maximal unramified extension of $F_0$ in $F_0^{\rm sep}$. If $F_1$ is
the maximal unramified extension of $F_0$ contained in $F$, we obtain a
decomposition

\begin{eq}\label{7.10}
F_1\otimes_{F_0}\Fsmile_0= \Fsmile_0\oplus\ldots\oplus \Fsmile_0\ \ ,
\end{eq}
corresponding to the $f=[F_1:F_0]$ different embeddings $\alpha: F_1\to
\Fsmile_0$ of $F_1$ into $\Fsmile_0$. For fixed $\alpha$, the extension
$\Fsmile_{\alpha}=F\otimes_{F_{1,\alpha}}\Fsmile_0$ is a totally ramified
field extension of degree $e=[F:F_0]/f$ of $\Fsmile_0$. The simple central
algebra $B\otimes_F\Fsmile_{\alpha}=B\otimes_{F_{1,\alpha}}\Fsmile_0$
splits, i.e.\ is isomorphic to a matrix algebra over $\Fsmile_{\alpha}$.
Similarly, we obtain a decomposition

\begin{eq}\label{7.11}
V\otimes_{F_0}\Fsmile_0= \bigoplus_{\alpha}V_{\alpha}\ \ ,
\end{eq}
where $V_{\alpha}$ is a $\Fsmile_{\alpha}$-vector space, all of the same
dimension $d$.

Since $F_1$ is an unramified extension we obtain similar decompositions
for the rings of integers, e.g.\

\begin{eq}\label{7.12}
\O_B\otimes_{\O_{F_0}}\O_{\Fsmile_0}=\bigoplus_{\alpha}\O_B\otimes_{\O_{F_1},\alpha}\O_{\Fsmile_0}\
\ ,
\end{eq}

\noindent where for each $\alpha$ the summand
$\O_B\otimes_{\O_{F_1},\alpha}\O_{\Fsmile_0}=
\O_B\otimes_{\O_F}\O_{\Fsmile_{\alpha}}$ is a parahoric order in the
matrix algebra $B\otimes_F\Fsmile_{\alpha}$. Similarly, the periodic
lattice chain ${\cal L}$ corresponds to periodic
$\O_{\Fsmile_{\alpha}}$-lattice chains ${\cal L}_{\alpha}$ in each
$\Fsmile_{\alpha}$-vector space $V_{\alpha}$ in (\ref{7.11}). Let $K$ be a
sufficiently big Galois extension of $F_0$. Then

\begin{eq}\label{7.13}
V\otimes_{F_0}K= \bigoplus_{\alpha}V_{\alpha}\otimes_{\Fsmile_0}\Ksmile\ \
,
\end{eq}
where $\Ksmile =\Fsmile_0.K$, and for each $\alpha$

\begin{eq}\label{7.14}
V_{\alpha}\otimes_{\Fsmile_0}\Ksmile =V_{\alpha,0}\oplus V_{\alpha,1}\ \ .
\end{eq}

Let $\Esmile =E.\Fsmile_0$. It now follows that our functor is
representable by a projective scheme over $\O_E$ which after base change
from $\O_E$ to $\O_{\Fsmile}$ becomes isomorphic to the product over all
$\alpha$ of schemes considered before for the data
$(\Fsmile_{\alpha}/\Fsmile_0, V_{\alpha}, \mu_{\alpha}, {\cal
L}_{\alpha})$.

Let us denote by $M^{\rm naive}$ the scheme over $\Spec\O_E$ associated in
this way to the data fixed in the beginning of this section. This is what
we call {\it a naive local model of $EL$-type.} We know from the special
case of a standard model for $GL_d$ that $M^{\rm naive}$ is rarely flat
over ${\rm Spec}\ \O_E$. We now define a closed subscheme $M^{\rm loc}$ of
$M^{\rm naive}$ which stands a better chance of being flat over the base
scheme.

Assume first that $B=F$ is a totally ramified extension. In the notation
introduced after (\ref{7.5}) let ${\cal L}={\cal L}_I$. For every $i\in I$
we obtain a morphism

\begin{eq}\label{7.15}
\pi_i: M^{\rm naive}\longrightarrow M^{\rm naive}(\Lambda_i)\ \ i\in I\ \
.
\end{eq}
Here $M^{\rm naive}(\Lambda_i)=M(\Lambda_i, \mu)$ denotes the standard
model associated to $\Lambda_i$ (and $(F,V,\mu)$). We then define in this
case

\begin{eq}\label{7.16}
M^{\rm loc}=\bigcap_{i\in I} \pi_i^{-1}(M^{\rm loc}(\Lambda_i))
\end{eq}
(scheme-theoretic intersection inside $M^{\rm naive}$).

In the general case we have, with the notation used above,

\begin{eq}\label{7.17}
M^{\rm naive}\otimes_{\O_E}\O_{\Esmile}=\prod_{\alpha} M^{\rm naive}
(\Fsmile_{\alpha} /\Fsmile_0, V_{\alpha},\mu_{\alpha}, {\cal L}_{\alpha})\
\ .
\end{eq}
Let

\begin{eq}\label{7.18}
\Msmile^{\rm loc}=\prod_{\alpha}M^{\rm loc} (\Fsmile_{\alpha}/\Fsmile_0,
V_{\alpha}, \mu_{\alpha}, {\cal L}_{\alpha})\ \ .
\end{eq}
The descent datum on $M^{\rm naive}\otimes_{\O_E}\O_{\Esmile}$ respects
the closed subscheme $\Msmile^{\rm loc}$ and hence defines a closed
subscheme $M^{\rm loc}$ of $M^{\rm naive}$.

We conjecture that $M^{\rm loc}$ is flat over ${\rm Spec}\ \O_E$. This
would constitute the analogue of the result of G\"ortz [G] which confirms
the conjecture in the case when $F/F_0$ is unramified.

\medskip
\smallskip
\bigskip\bigskip\bigskip

\leftline{George Pappas\hfill Michael Rapoport} \leftline{Dept. of
Mathematics\hfill Mathematisches Institut} \leftline{Michigan State
University\hfill der Universit\"at zu K\"oln} \leftline{E. Lansing\hfill
Weyertal 86 -- 90} \leftline{MI 48824-1027\hfill 50931 K\"oln}
\leftline{USA\hfill Germany}\leftline{email: pappas@math.msu.edu\hfill
email: rapoport@mi.uni-koeln.de}

\end{document}